\documentclass[10pt]{article}

\usepackage{amsmath}
\usepackage{amsfonts}
\usepackage{graphicx}
\usepackage{hyperref}

\title{Unstructured Space-Time Finite Element Methods for Optimal Sparse
  Control of\\ Parabolic Equations\footnote{This work has been supported by Johann Radon Institute for
  Computational and Applied Mathematics (RICAM) during the special semester on
  optimization that took place between 14th October and 11th December,
  2019, at RICAM in Linz, Austria.}}
\author{
  Ulrich Langer\footnote{RICAM, Austrian Academy of Sciences, Altenberger
    Stra{\ss}e 69, 4040 Linz, Austria, email: ulrich.langer@ricam.oeaw.ac.at}
  \,\,\,
  Olaf Steinbach\footnote{Institut f\"{u}r Angewandte Mathematik, Technische
    Universit\"{a}t Graz, Steyrergasse 30, 8010 Graz, Austria, email: o.steinbach@tugraz.at}
  \,\,\,
  Fredi Tr\"{o}ltzsch\footnote{Institut f\"{u}r Mathematik, Technische
    Universit\"{a}t Berlin, Stra{\ss}e des 17. Juni 136, 10623 Berlin,
    Germany, email: troeltzsch@math.tu-berlin.de}
  \,\,\,
  Huidong Yang\footnote{RICAM, Austrian Academy of Sciences, Altenberger
    Stra{\ss}e 69, 4040 Linz, Austria, email: huidong.yang@ricam.oeaw.ac.at}
}
\date{\today}

\begin{document}
\maketitle

\begin{abstract}
We consider a space-time finite element method on fully unstructured
simplicial meshes for optimal sparse control of semilinear parabolic
equations. The objective is a combination of a standard quadratic
tracking-type functional including a Tikhonov regularization term and of the
$L^1$-norm of the control that accounts for its spatio-temporal sparsity. We
use a space-time Petrov-Galerkin finite element discretization for the
first-order necessary optimality system of the associated discrete optimal
sparse control problem. The discretization is based on a variational
formulation that employs piecewise linear finite elements simultaneously in
space and time. Finally, the discrete nonlinear optimality system that
consists of coupled forward-backward state and adjoint state equations is
solved by a semismooth Newton method.  
\end{abstract}
\textbf{Keywords:}{ space-time finite element method, optimal sparse control,
  semilinear parabolic equations} \\
\textbf{MSC 2010:}{ 49J20,  35K20, 65M60, 65M50, 65M15, 65Y05}

\section{Introduction}
Optimal sparse control with the $L^1$-norm of the control in the objective
functional and with linear elliptic state equations has been analyzed in \cite{Stadler}
about a decade ago. The method was extended to semilinear elliptic
optimal control problems in \cite{CasasHerzogWachsmuth12} and to problems
governed by elliptic equations with uncertain coefficients in
\cite{LiStadler}. In
\cite{CasasRyllTroltzsch,CasasRyllTroltzsch15,RyllLoberMartensEngelTroltzsch},
the authors investigated problems of optimal sparse control for the
Schl\"{o}gl and FitzHugh-Nagumo systems, where traveling wave fronts or spiral
waves were controlled. Directionally spatio-temporal optimal sparse control of
linear/semilinear elliptic/parabolic equations and its corresponding
optimality conditions  was studied in
\cite{CasasHerzogWachsmuth,HerzogStadlerWachsmuth}. Optimality conditions for
directionally sparse parabolic control problems without control constraints
were considered in \cite{CasasMateosRosch}; see also
\cite{CasasMateosRosch18}. The optimal sparse controls considered therein
exhibit sparsity in space, but not necessarily in time. We also mention
another class of optimal sparse control problems for parabolic equations with
controls  in measure spaces instead of $L^1$-spaces, see, e.g.,
\cite{BounlangerTrautmann,CasasClasonKunisch,CasasKunisch,CasasZuazua,KunischPieperVexler}. A
thorough review of the existing literature on this challenging topic is beyond
the scope of this work. Therefore, we refer to the recent survey article
\cite{Casas} on sparse solutions in optimal control of both elliptic and
parabolic equations and the references therein.  

Moreover, numerical approximations of optimal sparse controls of elliptic and
parabolic problems were of great interest. For example, the standard $5$-point
stencil was used in \cite{Stadler} for the discretization of elliptic optimal
sparse control problems. In \cite{CasasHerzogWachsmuth12}, the authors proved
rigorous error estimates for the finite element approximation of semilinear
elliptic sparse control problems with with box constraints on the control. In
the discrete optimality conditions, they used piecewise linear approximations
for the state and adjoint state, while a piecewise constant ansatz was applied
to the control and the subdifferential. We also mention the approximation of
sparse controls by piecewise linear functions in
\cite{CasasHerzogWachsmuth120}. Here, the authors adopt special quadrature
formulae to discretize the squared $L^2$-norm and the $L^1$-norm of the
control in the objective functional. This leads to an elementwise
representation of the control and of the subdifferential.  

Later, error estimates were derived for the space-time finite element
approximation of parabolic optimal sparse control problems without control
constraints in \cite{CasasMateosRosch}. The discretization was performed on
tensor-structured space-time meshes. In the associated discrete optimal
control problem, these authors used a spatio-temporal finite element ansatz
for the state  that consists of products of continuous, piecewise linear basis
functions in space, and piecewise constant basis functions in time. For the
control, they utilize space-time elementwise constant basis functions. The
same space as for the state was employed for the adjoint state in the discrete
optimality system. Improved approximation rates were achieved in
\cite{CasasMateosRosch18}. For the control discretization, they employ basis
functions that are continuous and piecewise linear in space and piecewise
constant in time. These methods can be re-interpreted as an implicit Euler
discretization of the spatially discretized optimality system. 

For optimal sparse control of the Schl\"{o}gl and FitzHugh-Nagumo models 
that were considered in \cite{CasasRyllTroltzsch}, 
a semi-implicit Euler-method in time and continuous piecewise linear finite 
elements in space were applied to both the state and adjoint state equations. 
Recently, in \cite{UzuncaKucukseyhanYucel}, for the optimal control of the
convective FitzHugh-Nagumo equation, the state and adjoint state equations were
discretized by a symmetric interior penalty Galerkin method in space 
and the backward Euler method in time. 

In contrast to the discretization methods discussed above, we apply
continuous space-time finite element approximations on fully unstructured
simplicial meshes for parabolic optimal sparse control problems with control
constraints. This can be seen as an extension of the Petrov-Galerkin
space-time finite element method proposed in \cite{Olaf} for parabolic
problems,  and in our recent work \cite{LangerOlafTroltzschYang} for parabolic
optimal control problems. This kind of unstructured space-time finite element
approaches has gained increasing interest; see, e.g.,
\cite{BankVassilevskiZikatanov,Behr,KaryofylliWendlingLoicMichelNorbertBehr,LangerNeumullerSchafelner,
SteinbachYang01,Toulopoulos,DanwitzKaryofylliVioletaHostersBehr}, 
and the survey article \cite{SteinbachYang02}. 

In comparison to the more conventional time-stepping methods or
tensor-structured space-time methods \cite{Gander,GanderNeumuller,NagelLogashenkoSchroderYang},
this unstructured space-time approach provides us with more flexibility in
constructing parallel space-time solvers such as parallel space-time algebraic
multigrid preconditioners \cite{LangerNeumullerSchafelner} or space-time
balancing domain decomposition by constraints (BDDC) preconditioners
\cite{LangerYang}. Moreover, it becomes more convenient to realize simultaneous space-time adaptivity on
unstructured space-time meshes
\cite{LangerNeumullerSchafelner,LangerOlafTroltzschYang,SteinbachYang} than
the other methods. Here, time is just considered as another spatial
coordinate. For more comparisons of our space-time finite element methods with
others, we refer to \cite{SteinbachYang02}.    

The remainder of this paper will be structured as follows: 
Section~\ref{sec:modelprom} describes a model optimal sparse control problem that
we aim to solve. 
Some preliminary existing results concerning optimality
conditions are given in Section~\ref{sec:prelres}. 
The space-time finite
element discretization of the associated discrete optimal control problem, the
resulting discrete optimality conditions, as well as the application of the
semismooth Newton iteration are discussed in Section~\ref{sec:femdis}. 
The applicability of our proposed method is confirmed by two
numerical examples in Section~\ref{sec:numexam}. 
Finally, some conclusions are drawn in Section~\ref{sec:con}.   

\section{The optimal sparse control model problem}\label{sec:modelprom}
We consider the optimal sparse control problem
\begin{equation}\label{eq:model}
\displaystyle\min_{z\in Z_{ad}}{\mathcal J} (z) := 
\frac{1}{2} \, \| u_z  - u_Q \|_{L^2(Q)}^2 +
  \frac{\varrho}{2} \, \|z\|_{L^2(Q)}^2 + \mu \, \|z\|_{L^1(Q)},
\end{equation}
where the admissible set is
\begin{equation}\label{eq:adset}
  Z_{ad}=\left\{ z\in L^{\infty} (Q) : a\leq z(x,t) \leq b \textup { for a.a. } (x,t)\in Q\right\},
\end{equation}
and $u_z$ is the unique solution of the state equation
\begin{equation}\label{eq:state}
  \begin{aligned}
    \partial_t u  - \Delta_x u   + R(u) & =  z  && \text{ in } Q:=\Omega\times(0,T), \\
    u &  =     0 && \text{ on }  \Sigma:=\partial\Omega\times(0,T),\\
    u  & =    u_0  && \text{ on }  \Sigma_0:=\Omega\times\{0\}.
  \end{aligned}
\end{equation}
Here, the spatial computational domain $\Omega\subset{\mathbb R}^d$, $d\in
\{1,2,3\}$, is supposed to be bounded and Lipschitz, $T>0$ is the fixed
terminal time, $\partial_t$ denotes the partial time derivative,
$\Delta_x=\sum_{i=1}^d\partial^2_{x_i}$ the spatial Laplacian, and the source
term $z$ acts as a distributed control in $Q$. Moreover, $u_Q\in L^2(Q)$ is a
given desired state. We further assume  
\begin{equation*}
  -\infty < a < 0 < b<+\infty,\;
  \varrho > 0,\;
  \mu > 0.
\end{equation*}
The nonlinear reaction term $R$ is defined by
\begin{equation*}
  R(u)=(u-u_1)(u-u_2)(u-u_3)
\end{equation*}
with given real numbers $u_1\leq u_2\leq u_3$. The functional $g:L^1(Q)\rightarrow {\mathbb R}$
defined by $g(\cdot)=\|\cdot\|_{L^1(Q)}$ is Lipschitz continuous
and convex\, but not Fr\'{e}chet differentiable. 
We notice that similar model problems for semilinear equations have been studied, 
e.g., in \cite{Casas,CasasHerzogWachsmuth}.  

\section{Preliminary results}\label{sec:prelres}
Let us recall some facts that are known from literature, e.g.,
\cite{Casas}. For all $z\in L^p(Q)$, $p>d/2+1$, the state equation
(\ref{eq:state}) has a unique solution $u_z\in W(0,T)\cap L^{\infty}(Q)$,
where 
\begin{equation}\label{eq:wspace}
  W(0,T)=\left\{v\in L^2(0,T; H_0^1(\Omega)) : \partial_t v\in L^2(0,T;H^{-1}(\Omega))\right\}.
\end{equation}
Here, $H_0^1(\Omega):=\{v\in
H^1(\Omega):v=0\textup{ on }\partial\Omega\}$. The mapping $z\mapsto u_z$ is
continuously Fr\'{e}chet differentiable in these spaces. The optimal control
problem has at least one optimal control that is denoted by $\bar{z}$; the
associated optimal state is denoted by $\bar{u}$.     

If $\bar{z}$ is a locally optimal control of
the model problem, then there exist a unique adjoint state $\bar{p}\in W(0,T)$
and some $\bar{\lambda}\in \partial g 
(\bar{u})\subset L^{\infty}(Q)$ such that $(\bar{u}, \bar{p}, \bar{z}, \bar{\lambda})$ solves 
the optimality system
\begin{subequations}\begin{alignat}{2}
    &\partial_t {u}  - \Delta_x {u}   + R({u})  =  {z} \text{ in } Q,\;{u}  = 0
    \text{ on }  \Sigma,\; {u}  =    u_0 \text{ on }
    \Sigma_0,  \\ 
    &-\partial_t {p} -\Delta_x {p} + R'(u){p}={u}-u_Q \text{ in } Q,\;{p}  = 0
    \text{ on }  \Sigma,\; {p}  =    0 \text{ on }  \Sigma_T, \\
    &\int_Q ({p}+\varrho {z} + \mu {\lambda})(v-{z}) \, dx \, dt \geq
    0\textup{ for all } v \in { Z_{ad}}, \label{eq:varinq}
 \end{alignat}
\end{subequations}
where $\Sigma_T:=\Omega\times\{T\}$. A detailed discussion of this optimality
system leads to the relations
\begin{subequations}\begin{eqnarray}
    &&\bar{z} (x,t)  =  0 \Leftrightarrow |\bar{p} (x,t)|\leq \mu, \\
    &&\bar{z} (x,t)  = \textup{\bf Proj}_{[a,b]}\left(-\frac{1}{\varrho}\left(
    \bar{p} (x,t) + \mu \bar{\lambda} (x,t)\right)\right),\\
    &&\bar{\lambda} (x,t)=\textup{\bf Proj}_{[-1,1]}\left(-\frac{1}{\mu}\bar{p} (x,t)\right),
\end{eqnarray}\end{subequations}
that hold for a.a. $(x,t)\in Q$; see,
e.g., \cite{CasasRyllTroltzsch}. Here, the projection $\textup{\bf
  Proj}_{[\alpha,\beta]}:{\mathbb R}\rightarrow [\alpha,\beta]$ is
defined by 
$\textup{\bf Proj}_{[\alpha,\beta]} (q)=\max\{\alpha,\min\{q,\beta\}\}$; see, e.g.,
\cite{Troltzsch}. The subdifferential of the $L^1${-}norm of the control is given as
follows:
\begin{equation*}
  \bar{\lambda} \in \partial g(\bar{z}) 
  \Leftrightarrow 
  \left\{
  \begin{array}{ll}
    \bar{\lambda} = 1  &  \textup{ if }\bar{z}(x,t)>0,\\
    \bar{\lambda} \in [-1, 1] \quad & \textup{ if } \bar{z}(x,t)=0,\\
    \bar{\lambda} = -1 &  \textup{ if }\bar{z}(x,t)<0
  \end{array}
  \right.
\end{equation*}
for a.a. $(x, t) \in Q$ and $\bar{z}\in L^{\infty}(Q)$. By these
relations, we obtain the following form of an optimal control: 
\begin{equation}
  \bar{z}=
  \begin{cases}
    a &\textup{ on }{\mathcal A}_a:=\{(x,t)\in Q : 
-\bar{p}(x,t)+\mu < \varrho a\},\\
    b &\textup{ on }{\mathcal A}_b:=\{(x,t)\in Q : -\bar{p}(x,t) - \mu > \varrho b\},\\
    0 &\textup{ on }{\mathcal A}_0:=\{(x,t)\in Q : |\bar{p}(x,t)|\leq\mu \},\\
    -\frac{1}{\varrho}(\bar{p}-\mu)&\textup{ on } {\mathcal I}_{-}:=\{(x,t)\in Q : \varrho a \leq  -\bar{p}(x,t)+\mu < 0 \},\\
    -\frac{1}{\varrho}(\bar{p}+\mu)&\textup{ on } {\mathcal I}_{+}:=\{(x,t)\in Q : 0 < -\bar{p}(x,t) - \mu \leq \varrho b \}.
  \end{cases}
\end{equation}
The set ${\mathcal A}_0$ accounts for  the sparsity of the control.  

Eliminating the control from the optimality system and using the projection
formulae above, we obtain the following system for the state and the adjoint
state: 
\begin{subequations}\label{eq:cpsystem}\begin{alignat}{2}
    & \partial_t {u}  - \Delta_x {u}   + R({u})  =
    \textup{\bf Proj}_{[a,b]}\left(-\frac{1}{\varrho}\left( {p}  + \mu
    \textup{\bf Proj}_{[-1,1]}\left(-\frac{1}{\mu}{p} \right) \right)\right) \text{ in } Q, \\ 
    & {u}  = 0 \text{ on }  \Sigma,\; {u}  =    u_0 \text{ on }
    \Sigma_0, \nonumber \\ 
    & -\partial_t {p} -\Delta_x {p} +
    {R{'}}({u}){p}={u}-u_Q \text{ in } Q,\\  
    & {p}  = 0
    \text{ on }  \Sigma,\; {p}  =    0 \text{ on }  \Sigma_T. \nonumber
  \end{alignat}
\end{subequations}
Let us define the Bochner spaces for the state and adjoint state variables as
follows:
\begin{equation*}
  \begin{aligned}
    &X_0:=L^2(0,T;H_0^1(\Omega))\cap H_{0,}^1(0,T;H^{-1}(\Omega))
    =\{v\in W(0, T), v=0\textup{ on } \Sigma_0\},\\
    &X_T:=L^2(0,T;H_0^1(\Omega))\cap H_{,0}^1(0,T;H^{-1}(\Omega))
    =\{v\in W(0,T), v=0\textup{ on } \Sigma_T\},\\
    &Y:=L^2(0,T;H_0^1(\Omega)).
  \end{aligned}
\end{equation*}
Our space-time variational formulation for the coupled system
(\ref{eq:cpsystem}) reads: Find ${u}\in X_0$ and ${p}\in X_T$ such that the
variational equations
\begin{subequations}\label{eq:weakcpsystem}\begin{eqnarray}
    &&\int_Q \partial_tu \, v \, dx \, dt + 
\int_Q\nabla_xu\cdot\nabla_x
v\,dx\, dt   + \int_Q R({u}) \, v\, dx \, dt \nonumber\\
    &&\quad\quad= \int_Q {\textup{\bf Proj}_{[a,b]}\left(-\frac{1}{\varrho}\left( {p}  + \mu
      \textup{\bf
        Proj}_{[-1,1]}\left(-\frac{1}{\mu}{p}\right)\right)\right)} \, 
v\, dx \, dt \label{eq:varstate}\\  
    &&-\int_Q u \, q \, dx \, dt -\int_Q\partial_t p \, q\, dx \, dt +
\int_Q \nabla_xp\cdot\nabla_x q\, dx \, dt +
    \int_Q R'(u)p \, q\,dx \, dt \nonumber\\
    &&\quad\quad= -\int_Q u_Q \, q\,dx \, dt
\end{eqnarray}\end{subequations}
hold for all $v,q\in Y$. This system is solvable, because the optimal control
problem has at least one solution. Due to the projection formula
on the right hand side of (\ref{eq:varstate}), we need special care to
discretize the optimality system. In fact, following the discretization scheme proposed in
\cite{AradaCasasTroltzsch}, we go back to the variational inequality
(\ref{eq:varinq}) and use piecewise constant approximation for the control in
order to derive first order necessary optimality conditions for the associated
discrete optimal control problem. We will discuss this in the forthcoming Section.    

\section{Space-time finite element discretization}\label{sec:femdis}
For the space-time finite element approximation of the optimal control 
problem \eqref{eq:model}, we consider an admissible triangulation 
${\mathcal T}_h(Q)$ of the space-time domain $Q$ into shape regular
simplicial finite elements $\tau$. Here, the mesh size $h$ is defined by
$h=\max_{\tau \in {\mathcal T}_h } h_\tau$ with $h_\tau$ being the diameter of
the element $\tau$; see, e.g., \cite{Braess,Steinbach}. For simplicity, 
we assume $\Omega$ to be a polygonal spatial domain. Therefore, the 
triangulation exactly covers $Q = \Omega \times (0,T)$.

Let $S_h^1(Q)$ be the space of continuous and piecewise linear functions
that are defined with respect to the triangulation ${\mathcal T}_h(Q)$. The
discrete variational form of the state equation (\ref{eq:state}) reads as
follows: Find $u_h\in X_{0,h}=S_h^1(Q) \cap X_0$ such that 
\begin{equation}\label{eq:disweakstate}
\int_Q \partial_t u_h \, v_h \, dx \, dt + 
\int_Q \nabla_x u_h \cdot \nabla_x v_h \, dx \, dt + 
\int_Q R(u_h) \, v_h \, dx \, dt =
\int_Q z \, v_h \, dx \, dt 
\end{equation}
is satisfied for all $v_h\in X_{0,h}$. For approximating the control $z$, we
define the space    
\begin{equation*}  
Z_h = \Big \{ z_h\in L^{\infty}(Q) : 
z_h \textup{ is  constant in each } \tau\in  {\mathcal T}_h \Big \}.
\end{equation*}
An element $z_h \in Z_h$ can be represented in the form 
\begin{equation*}
z_h = \sum_{\tau \in {\mathcal T}_h} z_\tau {\mathcal X}_\tau,
\end{equation*}
with ${\mathcal X}_\tau$ being the characteristic function of $\tau$. 
Moreover, the set of discrete admissible controls is defined by
\begin{equation*} 
Z_{ad,h} = \Big \{ z_h \in Z_h : a \leq z_h|_\tau \leq b \; 
\textup { for all } \tau\in{\mathcal T}_h \Big  \} .
\end{equation*}
Now, we consider the discrete optimal control problem
\begin{equation*}
\min_{z_h\in Z_{ad,h}} {\mathcal J}_h (z_h) := 
\frac{1}{2} \, \| u_{z_h}  - u_Q \|_{L^2(Q)}^2 +
\frac{\varrho}{2}  \, \|z_h\|_{L^2(Q)}^2 + \mu \, \|z_h\|_{L^1(Q)},
\end{equation*}
where $u_{z_h}$ denotes the solution of the discrete variational problem
(\ref{eq:disweakstate}) with the discrete control $z_h$. We assume that the
discrete optimal control problem has at least one locally optimal  
control that is denoted by $\bar{z}_h$. The associated state is 
denoted by $\bar{u}_h$. If 
$\bar{z}_h=\sum_{\tau\in{\mathcal X}_h}\bar{z}_\tau {\mathcal X}_\tau$ 
is a locally optimal control of the discrete optimal control problem, 
then there exists a unique adjoint state 
$\bar{p}_h\in X_{T,h} = S_h^1(Q) \cap X_T$ and 
$\bar{\lambda}_h \in \partial g(\bar{z}_h)$ such that 
$(\bar{u}_h, \bar{p}_h, \bar{z}_h, \bar{\lambda}_h)$ solves the discrete 
optimality system   
\begin{subequations}\label{eq:disweakcpsystem}\begin{eqnarray}
\int_Q \partial_t u_h \, v_h \, dx \, dt + 
\int_Q \nabla_x u_h \cdot \nabla_x v_h \, dx \, dt + 
\int_Q R(u_h) \, v_h \, dx \, dt  \nonumber \hspace*{1cm} && \\
= \int_Q \bar{z}_h \, v_h \, dx \, dt, 
\textup{ for all } v_h \in X_{0,h}, \label{eq:disvarstate} \\  
-\int_Q u_h \, q_h \, dx \, dt - 
\int_Q \partial_t p_h \, q_h \, dx \, dt + 
\int_Q \nabla_x p_h \cdot \nabla_x q_h \, dx \, dt \hspace*{1cm} && \nonumber \\ 
+ \int_Q R'(u_h)p_h \, q_h \, dx \, dt = -\int_Q u_Q \, q_h \, dx \, dt, 
\textup{ for all } q_h\in X_{T,h}, && \label{eq:disvaradjoint} \\
\int_Q \Big( p_h + \varrho z_h + \mu \lambda_h \Big)\Big( v_h  - z_h \Big)
\, dx \, dt \, \ge \, 0, \textup { for all } v_h\in Z_{ad,h}. \label{eq:disvarinq}
&&\end{eqnarray}\end{subequations}
The existence of a solution to the discretized optimality system 
will not be discussed in this paper. We tacitly assume that a locally unique
solution exists. Note that $\bar{\lambda}_h \in \partial g(\bar{z}_h)$ is
equivalent to the form
\begin{equation*}
\bar{\lambda}_h = \sum_{\tau\in{\mathcal T}_h} 
\bar{\lambda}_\tau {\mathcal X}_\tau \quad \textup{ with } \quad
  \left\{
  \begin{array}{ll}
    \bar{\lambda}_\tau = 1  &  \textup{ if }\bar{z}_\tau >0,\\
    \bar{\lambda}_\tau \in [-1, 1] \quad & \textup{ if } \bar{z}_\tau=0,\\
    \bar{\lambda}_\tau = -1 &  \textup{ if }\bar{z}_\tau<0.
  \end{array}
  \right.
\end{equation*}
Then, the inequality (\ref{eq:disvarinq}) can be represented as follows:
\begin{equation*}
\sum_{\tau\in{\mathcal T}_h}\left( \int_\tau \bar{p}_h \, dx \, dt + 
|\tau| \Big( \varrho \bar{z}_\tau + \mu \bar{\lambda}_\tau \Big) \right) 
\Big( \bar{v}_\tau-\bar{z}_\tau \Big) \, \geq \, 0 
\textup{ for all } a\leq \bar{v}_\tau   \leq b,
\end{equation*} 
which is recasted in the equivalent form
\begin{equation*}
\left( \int_\tau \bar{p}_h \, dx \, dt + |\tau| \Big(
\varrho \bar{z}_\tau + \mu \bar{\lambda}_\tau \Big)\right) 
\Big( \bar{v}_\tau-\bar{z}_\tau \Big) \, \geq  \, 0, \; 
a\leq \bar{v}_\tau  \leq b,
\end{equation*} 
for all $\tau\in{\mathcal T}_h$. Therefore, we have the projection
representation formula \cite{Troltzsch}
\begin{equation}\label{eq:cellproj}
\bar{z}_\tau = \textup{\bf Proj}_{[a,b]}\left( -\frac{1}{\varrho} \left(
\frac{1}{|\tau|} \int_\tau \bar{p}_h \, dx \, dt + 
\mu \bar{\lambda}_\tau\right)\right)
\end{equation}  
for the optimal control on each element $\tau\in {\mathcal T}_h$. From this, we
have the following results:
\begin{equation*}
\begin{aligned}
&\bar{z}_\tau = 0 \; \Leftrightarrow \;
\frac{1}{|\tau|} \left| \int_\tau \bar{p}_h \, dx \, dt \right| \leq \mu, \\
&\bar{\lambda}_\tau = \textup{\bf Proj}_{[-1,1]}\left(
-\frac{1}{\mu|\tau|}\int_\tau \bar{p}_h \, dx \, dt \right).
\end{aligned}
\end{equation*}
The above results are obtained by the discretization approach
proposed in \cite{AradaCasasTroltzsch, CasasHerzogWachsmuth12} for the optimal
control of semilinear elliptic equations, and in \cite{CasasRyllTroltzsch} of
the Schl\"{o}gl and FitzHugh-Nagumo systems. In fact, by a close look to the
projection formula, we have the following form of the discrete optimal
control:      
\begin{equation}
  \bar{z}_\tau =
  \begin{cases}
    a &\textup{ on }{\mathcal A}_{a,{\mathcal T}_h}:=\{\tau\in {\mathcal T}_h : -\bar{p}_{\tau}+\mu < \varrho a\},\\
    b &\textup{ on }{\mathcal A}_{b,{\mathcal T}_h}:=\{\tau\in  {\mathcal T}_h : -\bar{p}_{\tau} - \mu > \varrho b\},\\
    0 &\textup{ on }{\mathcal A}_{0,{\mathcal T}_h}:=\{\tau\in {\mathcal T}_h : |\bar{p}_{\tau}|\leq\mu \},\\
    -\frac{1}{\varrho}(\bar{p}_\tau -\mu)&\textup{ on } {\mathcal I}_{-,{\mathcal T}_h} :=\{\tau\in {\mathcal T}_h : \varrho a \leq  -\bar{p}_{\tau}+\mu < 0 \},\\
    -\frac{1}{\varrho}(\bar{p}_\tau +\mu)&\textup{ on } {\mathcal I}_{+,{\mathcal T}_h} :=\{\tau\in {\mathcal T}_h : 0 < -\bar{p}_{\tau} - \mu \leq \varrho b \},\\
  \end{cases}
\end{equation}
where 
\[
\bar{p}_{\tau}=\frac{1}{|\tau|}\int_\tau \bar{p}_{h} \, dx \, dt .
\] 
Inserting 
(\ref{eq:cellproj}) into (\ref{eq:disvarstate}), we obtain an equivalent form 
of the discrete optimality system that consists of the state and adjoint state 
equations. Namely, find $\bar{u}_h\in X_{0,h}$ and $\bar{p}_h\in X_{T,h}$ such 
that $(\bar{u}_h,\bar{p}_h)$ solves the coupled system
\begin{subequations}\label{eq:reduceddisweakcpsystem}\begin{eqnarray}
\int_Q \partial_t u_h \, v_h \, dx \, dt +
\int_Q \nabla_x u_h \cdot \nabla_x v_h \, dx \, dt + 
\int_Q R(u_h) \, v_h \, dx \, dt \hspace*{1cm} & & \nonumber \\
-\sum_{\tau\in {\mathcal T}_h} \int_{\tau} \textup{\bf Proj}_{[a,b]} 
\left(-\frac{1}{\varrho} \left( 
\bar{p}_\tau
+ \mu \, \textup{\bf Proj}_{[-1,1]} \left(
-\frac{1}{\mu}
\bar{p}_\tau
\right)\right)\right) \, v_h \, dx \, dt \hspace*{1cm} & & \nonumber \\
=0, \textup{ for all }v_h\in X_{0,h}, && \\ 
-\int_Q u_h \, q_h \, dx \, dt -\int_Q \partial_t p_h \, q_h \, dx \, dt +
\int_Q \nabla_x p_h \cdot \nabla_x q_h \, dx \, dt \hspace*{1cm} \nonumber && \\
+ \int_Q R'(u_h)p_h \, q_h \,dx \, dt \, = \, - \int_Q u_Q \, q_h \, dx \, dt, 
\textup{ for all } q_h\in X_{T,h}.
\end{eqnarray}
\end{subequations}
The convergence of the solution of the discrete optimal control problem to the
solution of its associated continuous optimal control problem as well as the
error analysis of our finite element approximation are beyond the scope of
this work and will be studied elsewhere.   

To solve the above discrete coupled nonlinear optimality system,
we apply the semismooth Newton method as discussed in
\cite{Stadler}, where a generalized derivative needs to be computed at each
Newton iteration. In fact, each iteration turns out to be one step of a
primal-dual active set strategy \cite{ItoKunisch}: Given $(u_h^k, p_h^k)$,
find $(\delta u_h, \delta p_h)$ such that    
\begin{alignat*}{2} 
& \int_Q \partial_t \delta u_h \, v_h \,dx \, dt +
\int_Q \nabla_x \delta u_h \cdot \nabla_x v_h \, dx \, dt + 
\int_Q R'(u_h^k) \, \delta u_h \, v_h \, dx \, dt \\
& + \frac{1}{\varrho} \left(
{\mathcal X}_{{\mathcal I}_{-,{\mathcal T}_h}} + 
{\mathcal X}_{{\mathcal I}_{+,{\mathcal T}_h}} \right)
\left( \sum_{\tau\in {\mathcal T}_h} \int_\tau
\left(\frac{1}{|\tau|} \int_\tau \delta p_h \, dx \, dt \right) 
v_h \, dx \, dt \right)  \\ 
& = -\int_Q \partial_t u_h^k \, v_h \, dx \, dt -
\int_Q \nabla_x u_h^k \cdot \nabla_x v_h \, dx \, dt + 
\int_Q R(u_h^k) \, v_h \, dx \, dt \\
& + {\mathcal X}_{{\mathcal A}_{a,{\mathcal T}_h}} 
\left( \sum_{\tau\in {\mathcal T}_h} \int_\tau a \, v_h \, dx \, dt \right) +
{\mathcal X}_{{\mathcal A}_{b,{\mathcal T}_h}}
\left( \sum_{\tau\in {\mathcal T}_h} \int_\tau b \, v_h \, dx \, dt \right) \\ 
& -\frac{1}{\varrho}{\mathcal X}_{{\mathcal I}_{-,{\mathcal T}_h}}
\left( \sum_{\tau\in {\mathcal T}_h} \int_\tau \left(
\frac{1}{|\tau|} \int_\tau p_h^k \, dx \, dt - \mu \right) v_h \, dx \, dt 
\right) \\  
& -\frac{1}{\varrho}{\mathcal X}_{{\mathcal I}_{+,{\mathcal T}_h}}
\left( \sum_{\tau\in {\mathcal T}_h} \int_\tau
\left( \frac{1}{|\tau|} \int_\tau p_h^k \, dx \, dt + \mu \right)
v_h \, dx \, dt \right)
\end{alignat*}
and
\begin{alignat*}{2} 
&-\int_Q \delta u_h \, q_h \, dx \, dt -
\int_Q \partial_t \delta p_h \, q_h \, dx \, dt +
\int_Q \nabla_x \delta p_h \cdot \nabla_x q_h \, dx \, dt \\
& + \int_Q R'(u_h^k) \delta p_h \, q_h \, dx \, dt +
\int_Q R''(u_h^k)p_h^k \delta u_h \, q_h \, dx \, dt  \\
& =  \int_Q u_h^k \, q_h \, dx \, dt +
\int_Q \partial_t p_h^k \, q_h \, dx \, dt -
\int_Q \nabla_x p_h^k \cdot \nabla_x q_h \, dx \, dt \\
& - \int_Q R'(u_h^k)p_h^kq_h \, dx \, dt 
- \int_Q u_Q \, q_h \, dx \, dt
\end{alignat*}
are fulfilled, and $u_h^{k+1}=u_h^k+\omega\delta u_h$,
$p_h^{k+1}=p_h^k+\omega\delta p_h$ with some damping parameter
$\omega\in(0,1]$. 

\section{Numerical experiments}\label{sec:numexam}
For the two numerical examples considered in this section, we set
$\Omega=(0,1)^2$, $T=1$, and therefore $Q=(0,1)^3$. Using an 
octasection-based refinement \cite{Jurgen}, we uniformly decompose 
the space-time cylinder $Q$ until the mesh size reaches $h=1/128$.
Therefore, the total number of degrees of freedom for the
coupled state and adjoint state equations is $4,194,304$. 
We will also use an adaptive refinement procedure that is driven by a
residual based error indicator for the coupled state and adjoint state system 
similar to that one that was developed for the state equation in 
\cite{SteinbachYang}. We perform our numerical tests on a desktop with 
Intel@ Xeon@ Processor E5-1650 v4 ($15$ MB Cache, $3.60$ GHz), 
and $64$ GB memory. For the nonlinear first order necessary optimality 
system, we use the relative residual error $10^{-5}$ as a stopping criterion 
in the semismooth Newton iteration, 
whereas the algebraic multigrid preconditioned GMRES solver for the linearized
system at each Newton iteration is stopped after a residual error reduction by 
$10^{-6}$; cf. also \cite{SteinbachYang02}.

\subsection{Moving target (Example 1)}
\begin{figure}
  \centering
  \includegraphics[scale=0.185]{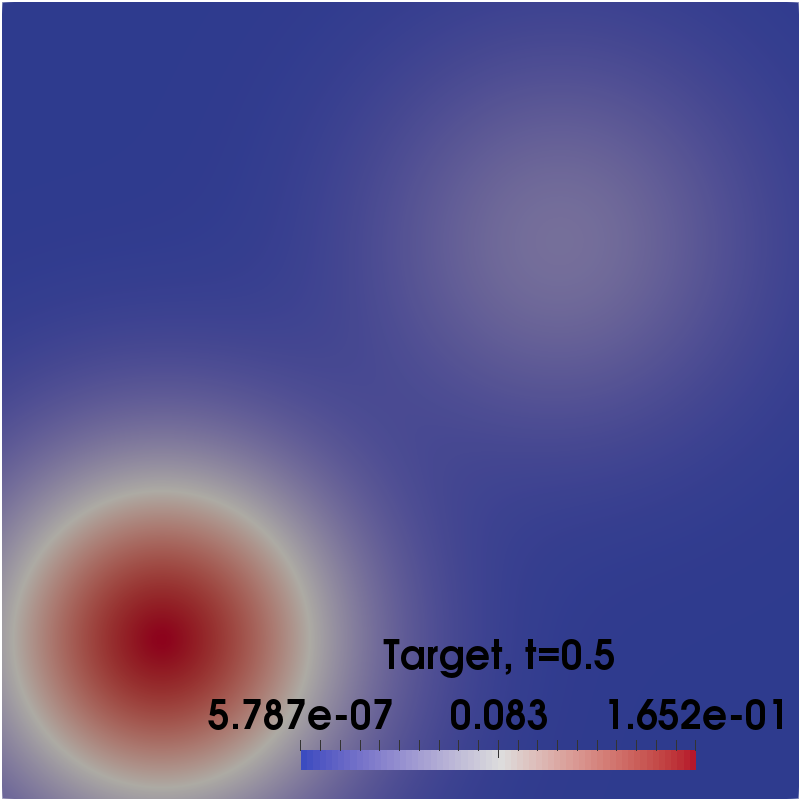}
  \includegraphics[scale=0.185]{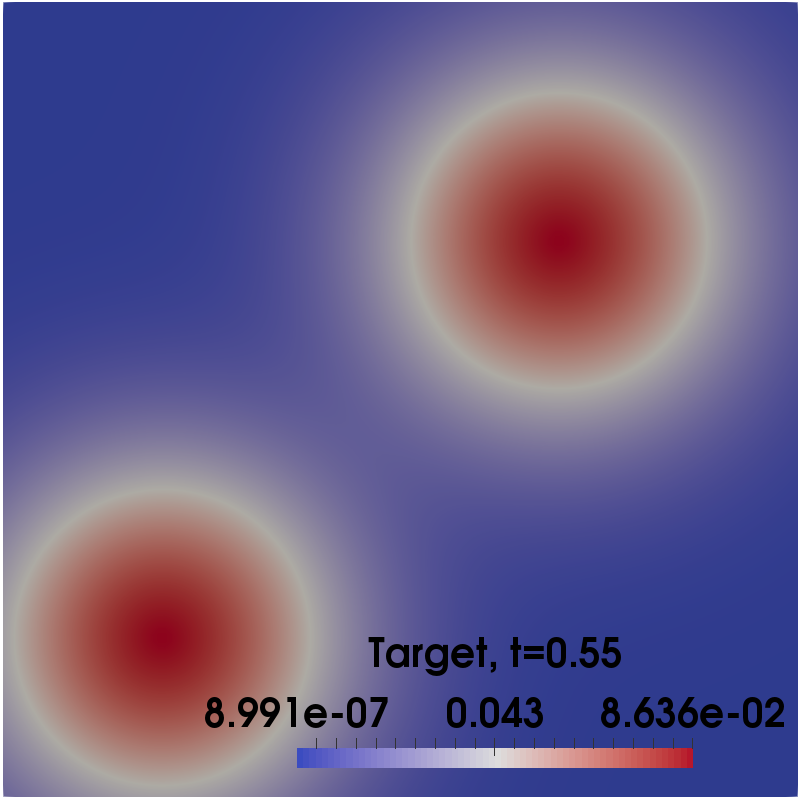}
  \includegraphics[scale=0.185]{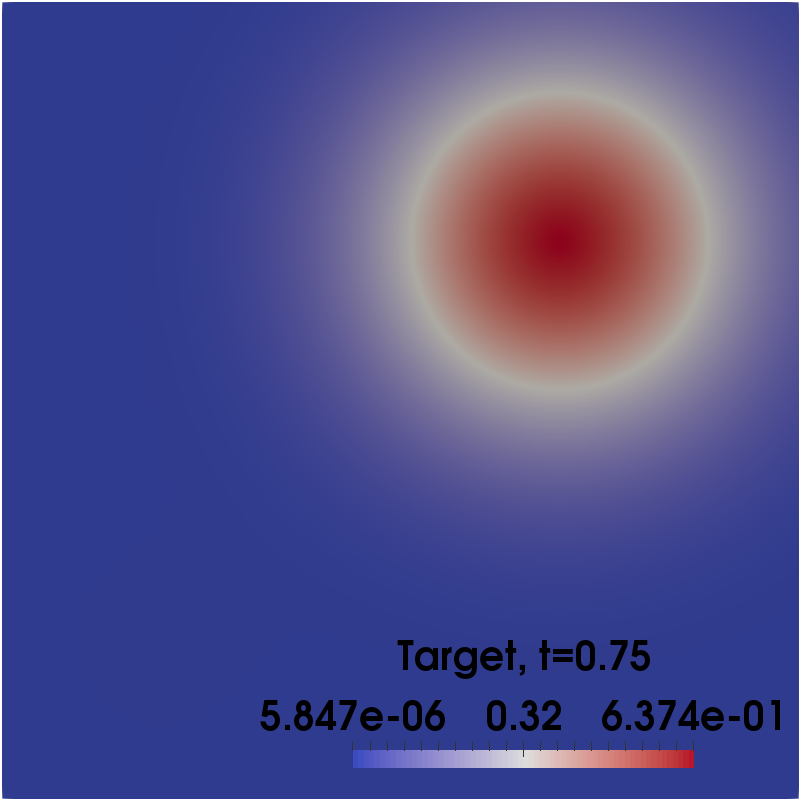}
  \caption{Example 1, plots of the target at time $t=0.5,\;0.55,\;0.75$
    for the moving target example.} \label{fig:targetexp}
\end{figure}

In the first example, the desired state is given by the function
\begin{equation*}
  \begin{aligned}
  u_Q(x,t)=&\exp\left(
  -20(x_1-0.2)^2+(x_2-0.2)^2+(t-0.2)^2
  \right)+\\
  &\quad\quad\exp\left(
  -20(x_1-0.7)^2+(x_2-0.7)^2+(t-0.9)^2
  \right),
  \end{aligned}
\end{equation*}
which is adapted from an example constructed in
\cite{CasasHerzogWachsmuth}; see an illustration of the target at time
$t=0.5$, $0.55$, and $0.75$ in Fig.~\ref{fig:targetexp}. The same desired state
was also used in the numerical test for spatially directional sparse control
in \cite{CasasMateosRosch,CasasMateosRosch18}. The parameters in the optimal
control problem are $\varrho=10^{-4}$, $\mu=0.004$, $a=-10$, and $b=20$. For the
nonlinear reaction term in the state equation, we set
$R(u)=u(u-0.25)(u+1)$. Homogeneous initial and Dirichlet boundary conditions
are used for the state equation. 

\begin{figure}
\centering
\includegraphics[scale=0.25]{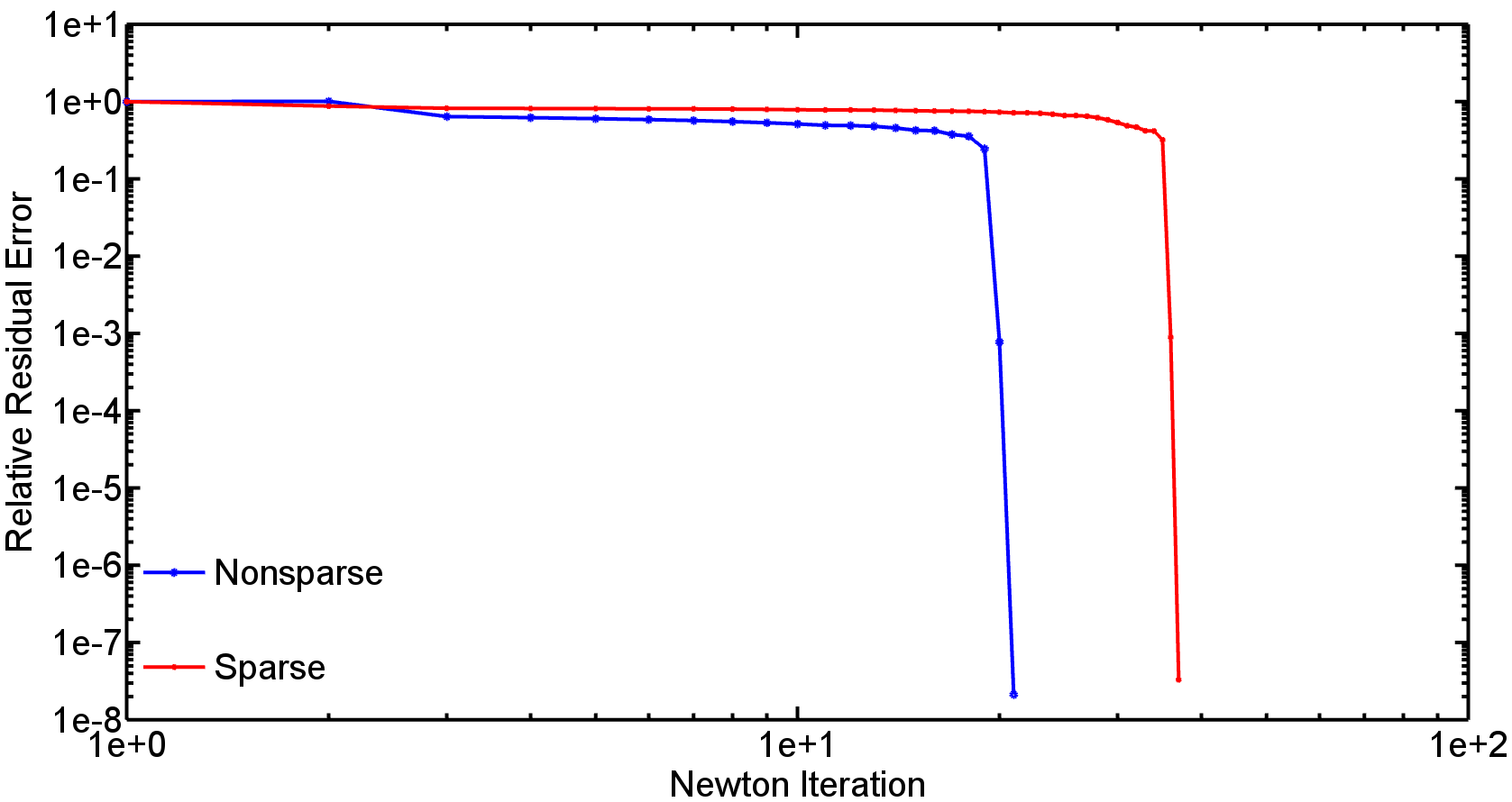}
\caption{Example 1, relative residual error reduction in the semismooth
Newton method for the nonsparse and sparse control in the moving 
target example.} \label{fig:expnewton}
\end{figure}

\noindent
To reach the relative residual error $10^{-5}$ for the nonlinear first order
necessary optimality system, we needed $21$ and $37$ semismooth Newton
iterations for the nonsparse and sparse optimal control, respectively; see
Fig.\ref{fig:expnewton}. We clearly see a superlinear convergence of the
semismooth Newton method. The total system assembling and solving time is about
$4$ and $4.7$ hours on the desktop computer, respectively.
 
\begin{figure}
  \centering
  \includegraphics[scale=0.18]{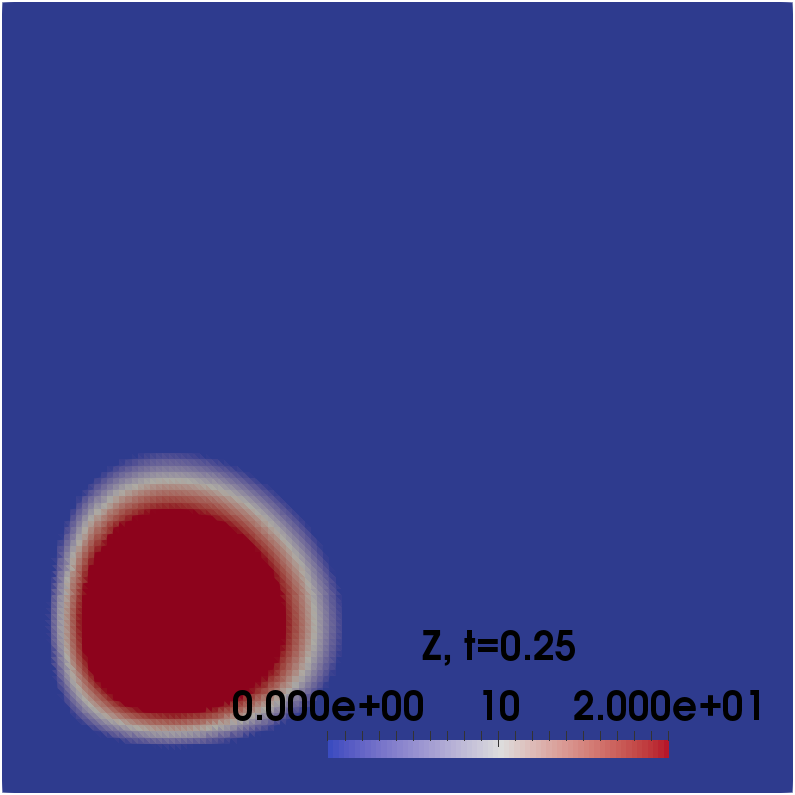}
  \includegraphics[scale=0.18]{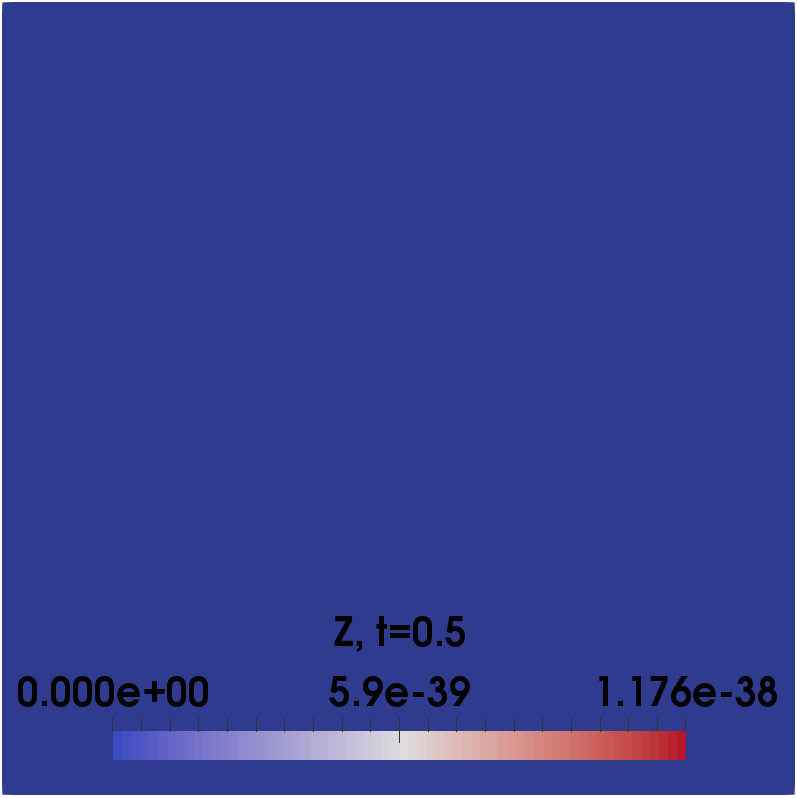}
  \includegraphics[scale=0.18]{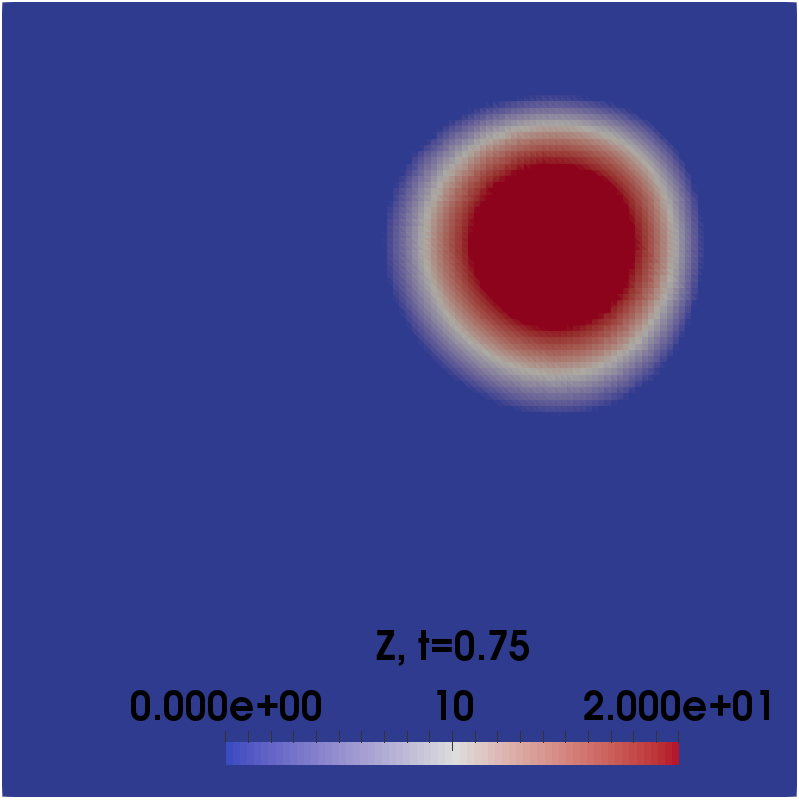}\\
  \includegraphics[scale=0.18]{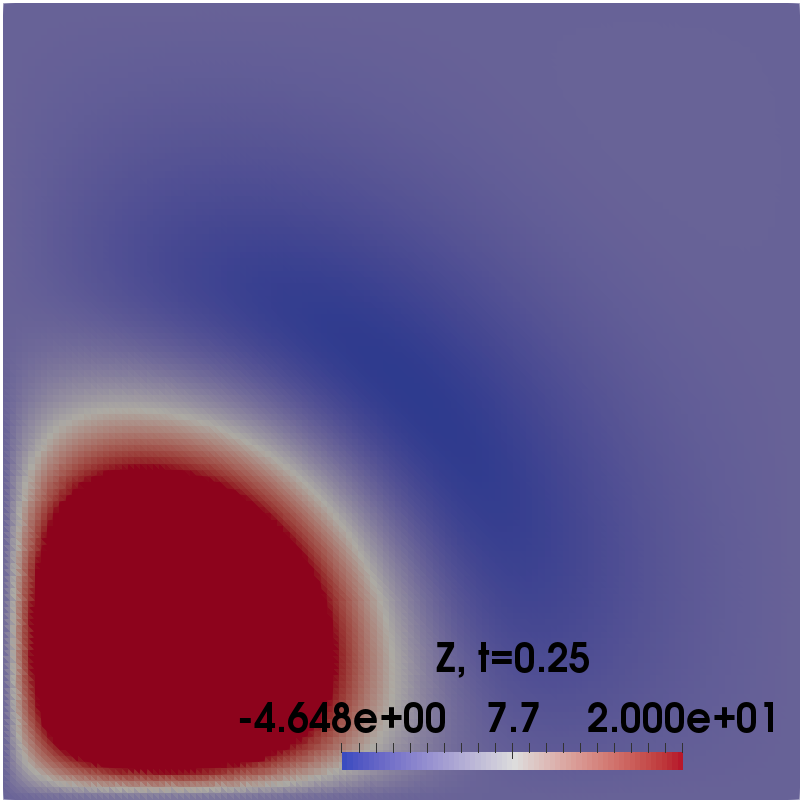}
  \includegraphics[scale=0.18]{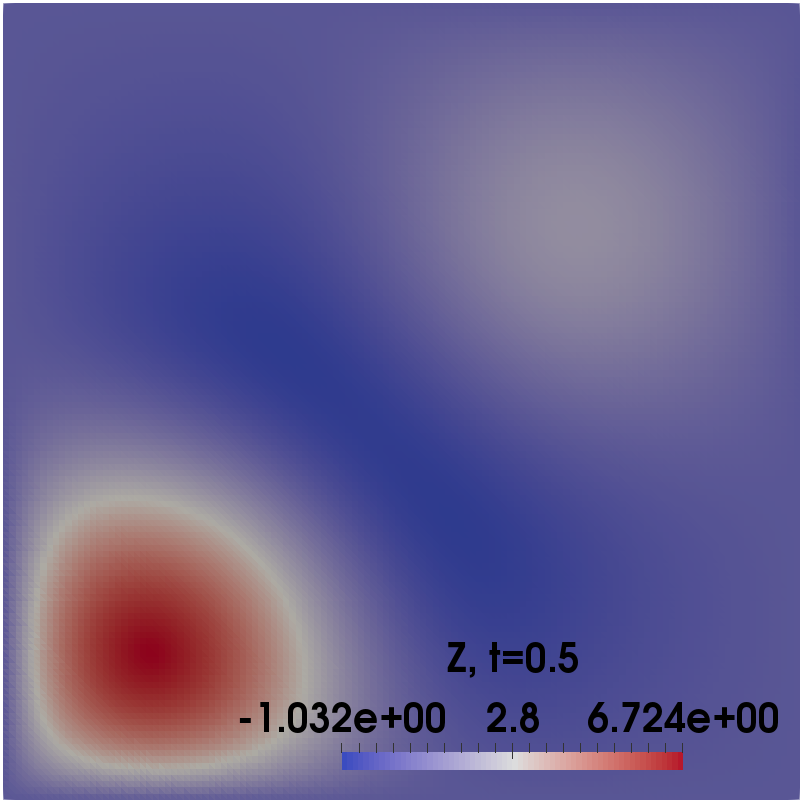}
  \includegraphics[scale=0.18]{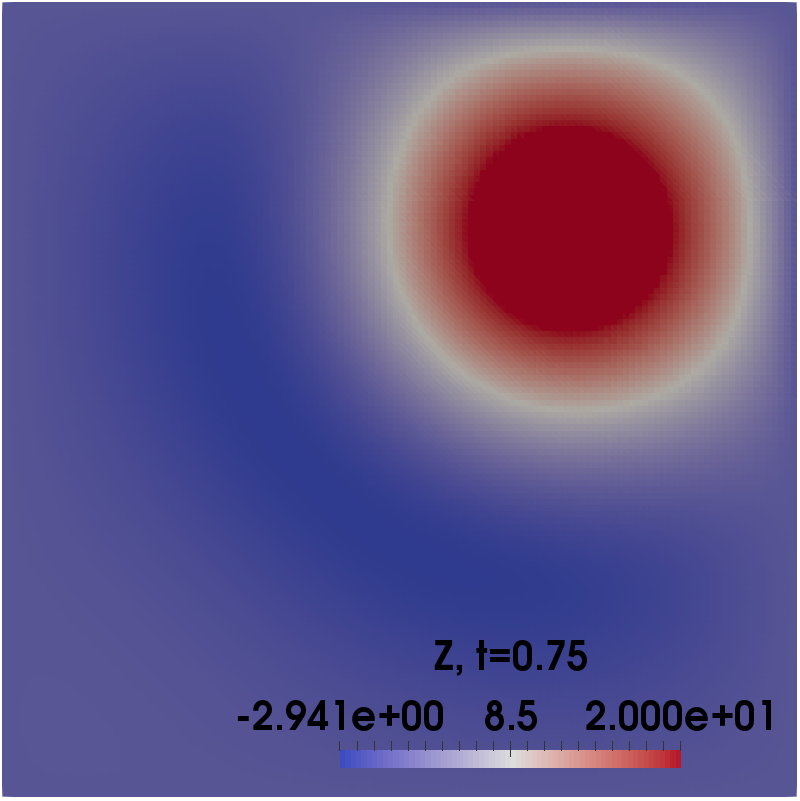}
  \caption{Example 1, Comparisons of sparse (up) and nonsparse controls (down) at time $t=0.25,\;0.5,\;0.75$
    for the moving target example.} \label{fig:controlexp}
\end{figure}

\begin{figure}
\centering
  \includegraphics[scale=0.25]{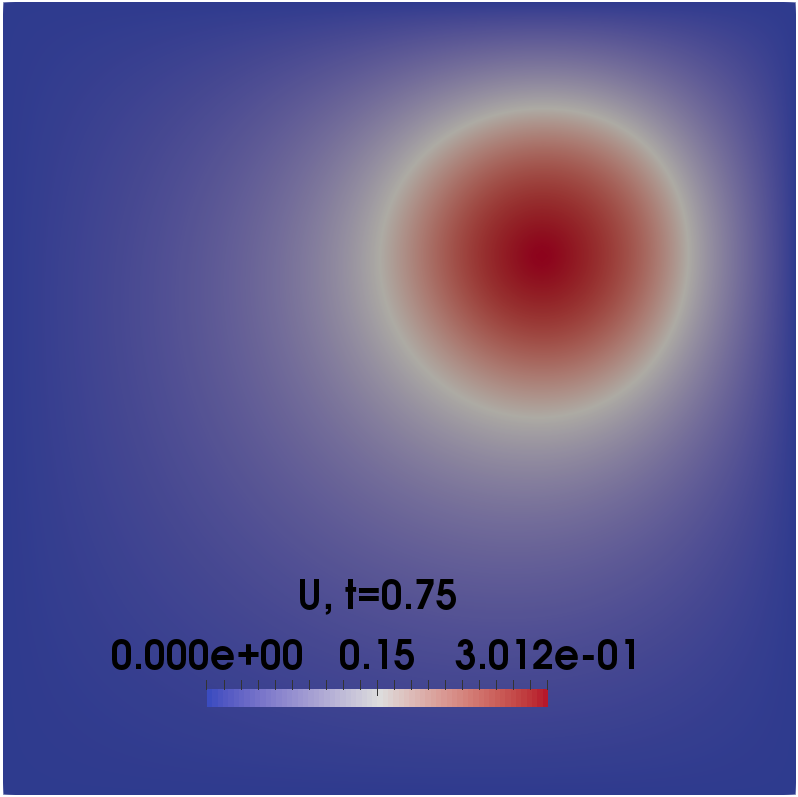}
  \includegraphics[scale=0.25]{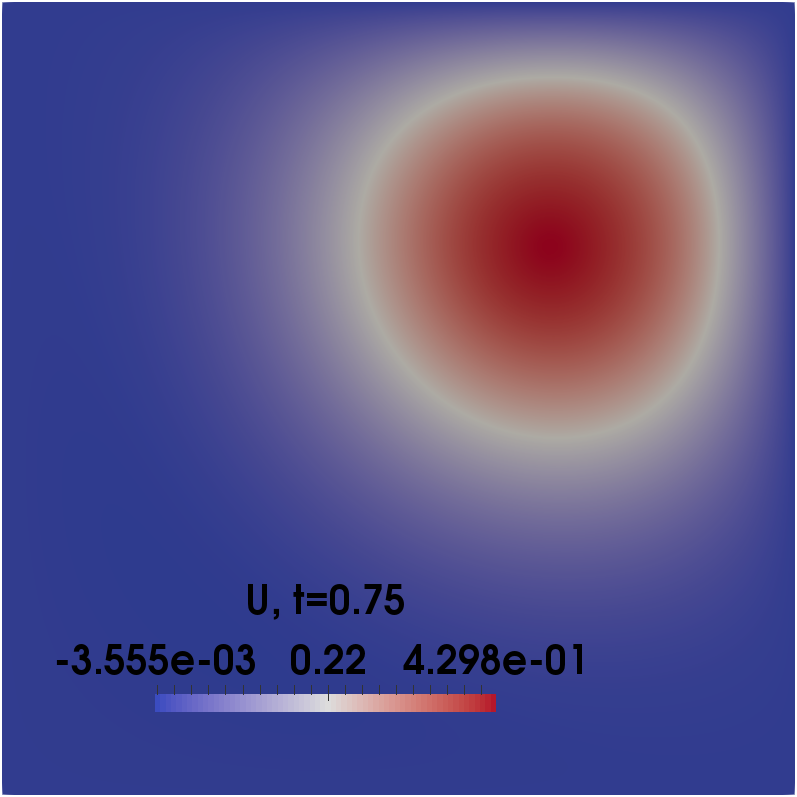}
  \caption{Example 1, Comparisons of states associated to sparse (left) and nonsparse controls (right) at time $t=0.75$ for the moving target example.} \label{fig:stateexp}
\end{figure}

\begin{figure}
  \centering
  \includegraphics[scale=0.26]{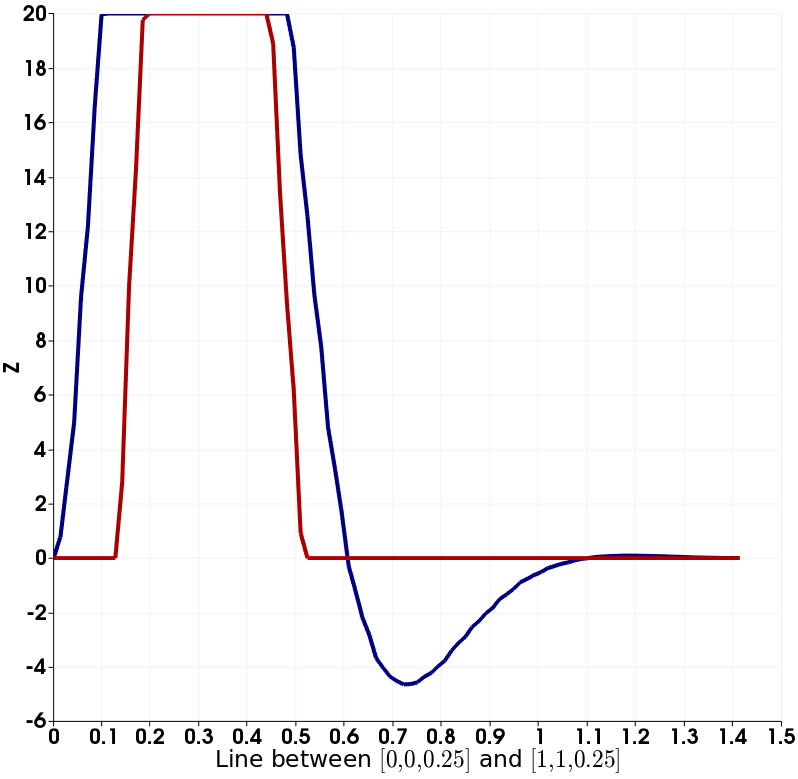}
  \includegraphics[scale=0.26]{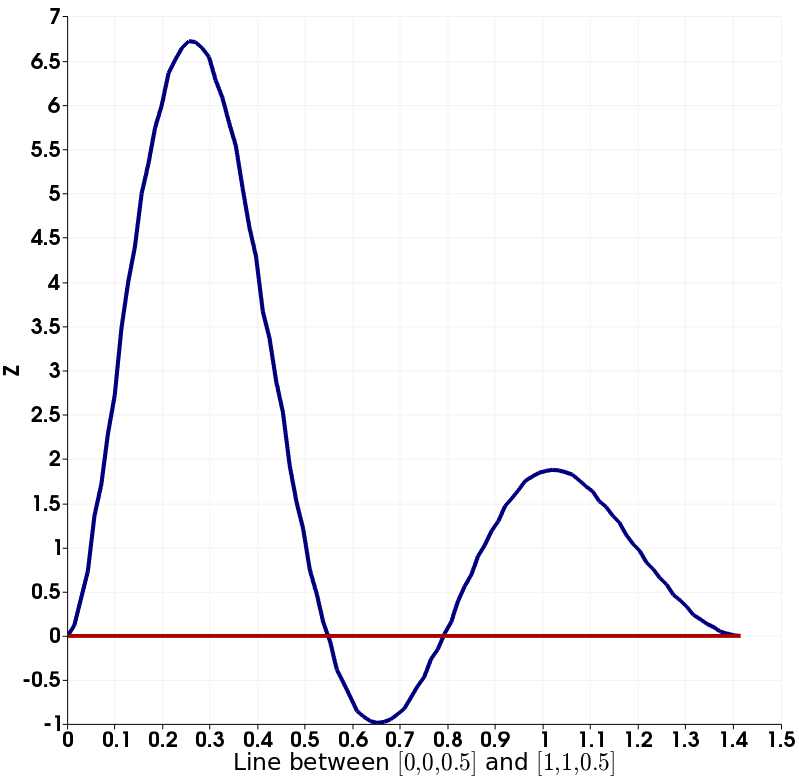}
  \includegraphics[scale=0.26]{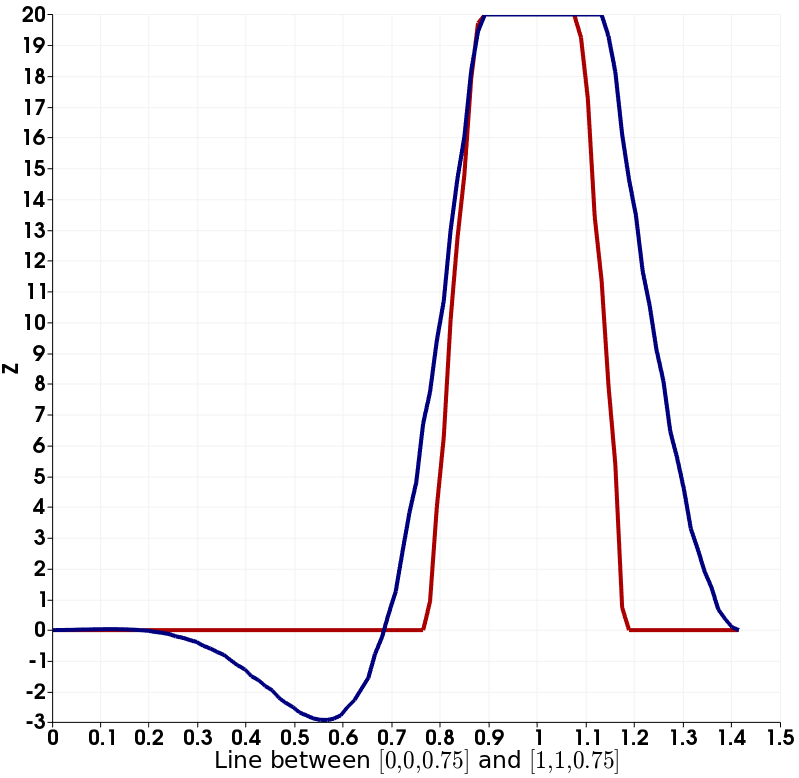}
  \caption{Example 1, comparisons of sparse (red) and nonsparse (blue)
    controls along the line between $[0, 0, 0.25]$ and $[1, 1, 0.25]$,
    between $[0, 0, 0.5]$ and $[1, 1, 0.5]$, and between $[0, 0, 0.75]$
    and $[1, 1, 0.75]$ for the moving target example.} \label{fig:expsparsecontrolline}  
\end{figure}

\noindent
Comparisons of sparse and nonsparse controls as well as of associated states at
different times are displayed in Fig.~\ref{fig:controlexp} and
Fig.~\ref{fig:stateexp}, respectively. A close look to sparse and nonsparse
controls along different lines in the space-time domain is illustrated in
Fig.~\ref{fig:expsparsecontrolline}. In this example, we clearly see that the
$L^1$ cost functional promotes spatial and temporal sparsity, with some
precision loss of the associated state to the 
target, cf. Fig.~\ref{fig:stateexp}.   

\subsection{Turning wave target (Example 2)}
\begin{figure}
  \centering
  \includegraphics[scale=0.2]{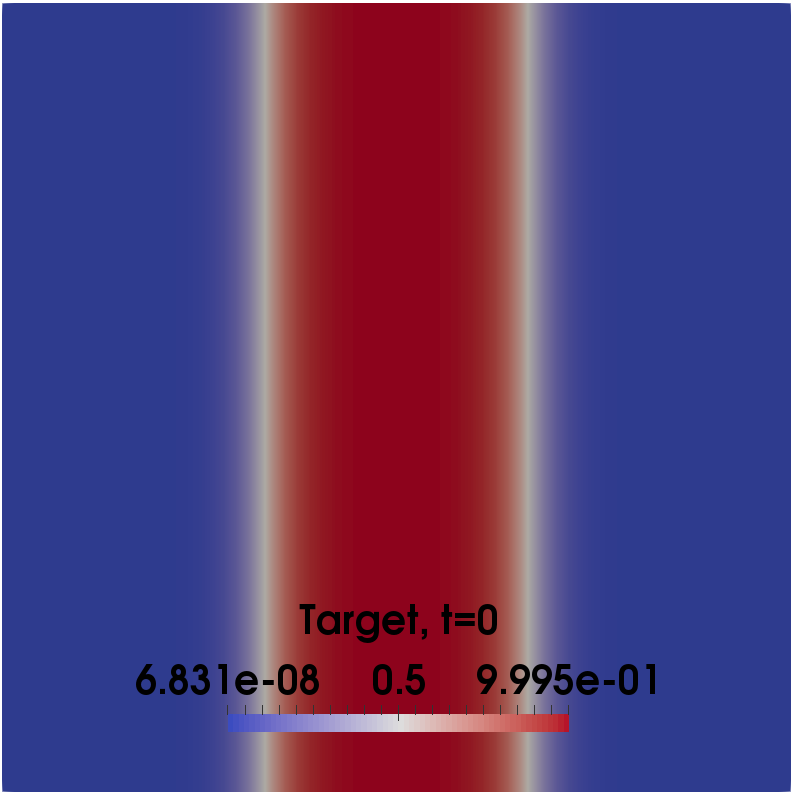}
  \includegraphics[scale=0.2]{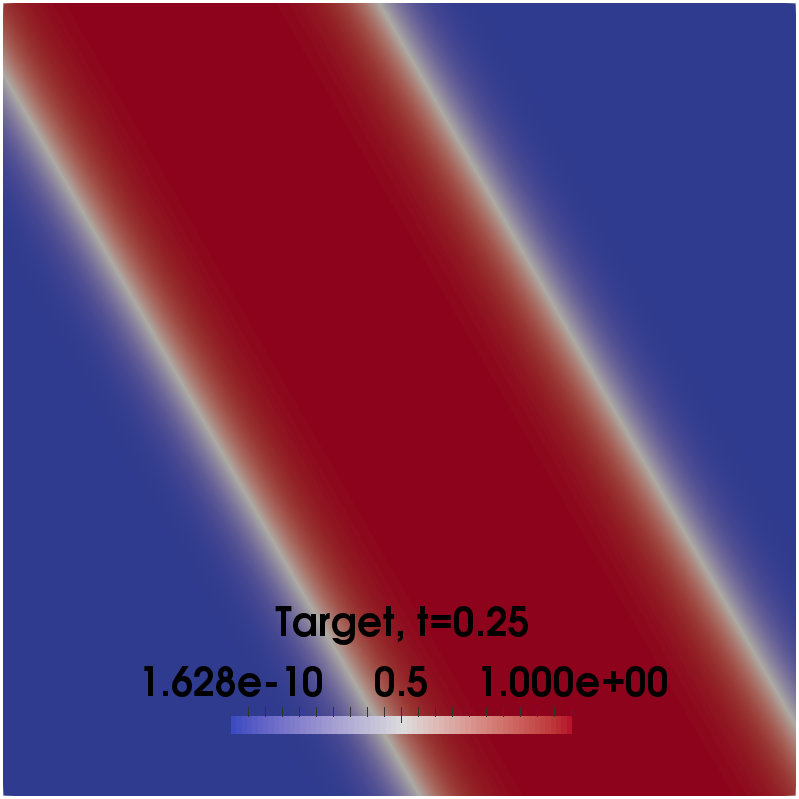}
  \includegraphics[scale=0.2]{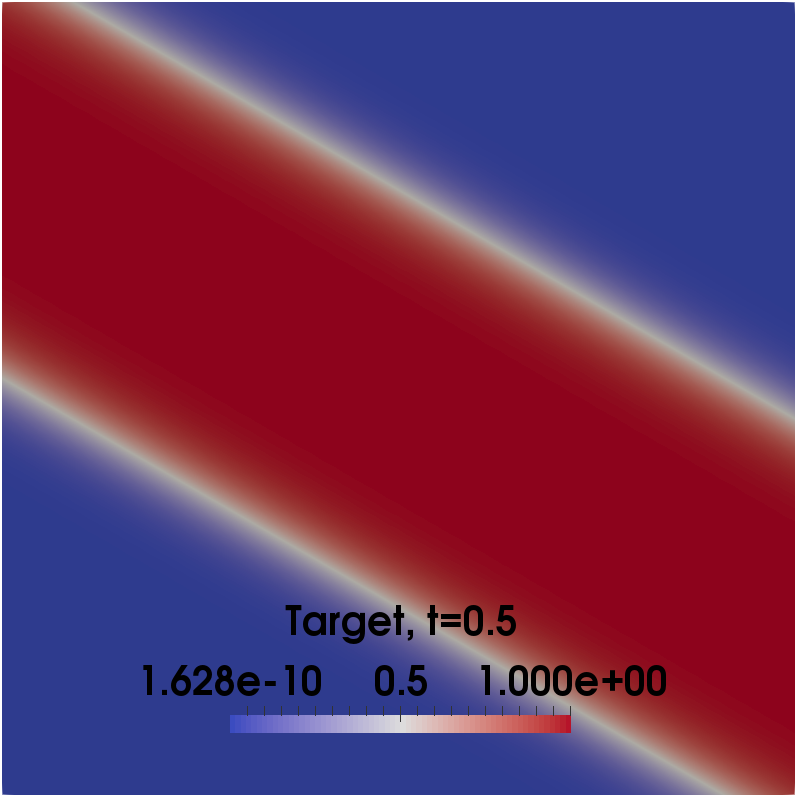}
  \includegraphics[scale=0.2]{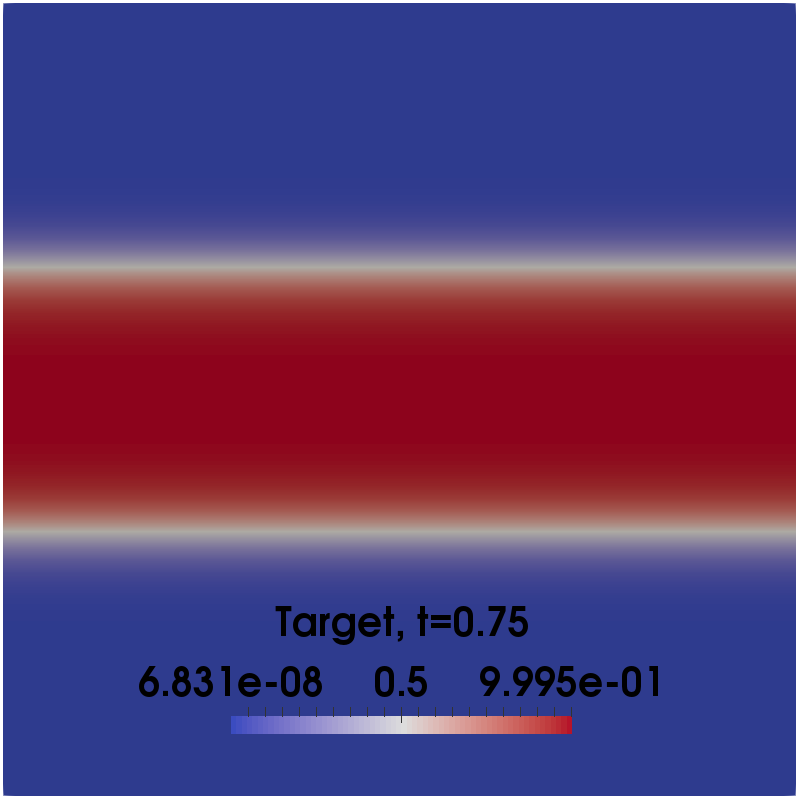}
  \caption{Example 2, plots of the target at time $t=0,\;0.25,\;\;0.5,\;0.75$
    for the turning wave example.} \label{fig:targetturningwave}
\end{figure}

In this example, we consider the target
\begin{equation*}
  \begin{aligned}
    &{u}_Q(x,t)=\left( 1.0+ \exp \left( \frac{\cos ( g(t)
      )\left(\frac{70}{3}- 70x_1\right)+ \sin ( g(t) ) \left(\frac{70}{3} -
      70x_2\right) }{\sqrt{2}} \right) \right)^{-1}\\
    &\quad + \left( 1.0+ \exp \left( \frac{ \cos ( g(t) ) \left( 70x_1 -
      \frac{140}{3}\right) + \sin ( g(t) ) \left( 70x_2 -
      \frac{140}{3}\right)}{\sqrt {2}} \right) \right)^{-1} - 1,
  \end{aligned}
\end{equation*}
where $g(t)=\frac{2\pi}{3}\min\left\{\frac{3}{4}, t\right\}$. This is an
adapted version of the turning wave example considered in
\cite{CasasRyllTroltzsch}.  The wave front turns $90$ degrees from time $t=0$
to $t=0.75$, and remains fixed after $t=0.75$; see the target at $t=0$, $0.25$,
$0.5$, and $0.75$ as illustrated in Fig.~\ref{fig:targetturningwave}. 

The nonlinear reaction term is given by $R(u)=u(u-0.25)(u+1)$. 
We use the initial data
\[
u_0(x)=\left(1+\exp\left(\frac{\frac{70}{3}-70x_1}{\sqrt{2}}\right)\right)^{-1}
+\left(1+\exp\left(\frac{70
x_1-\frac{140}{3}}{\sqrt{2}}\right)\right)^{-1}-1
\] 
on $\Sigma_0$, and homogeneous Neumann boundary condition on $\Sigma$ for 
the state. As parameters, we use $\varrho=10^{-6}$, $\mu=10^{-4}$ for the sparse 
case and $\varrho=10^{-6}$, $\mu=0$ for the nonsparse case. The bounds $a=-100$ 
and $b=100$ are set for both cases. 

\begin{figure}
\centering
\includegraphics[scale=0.25]{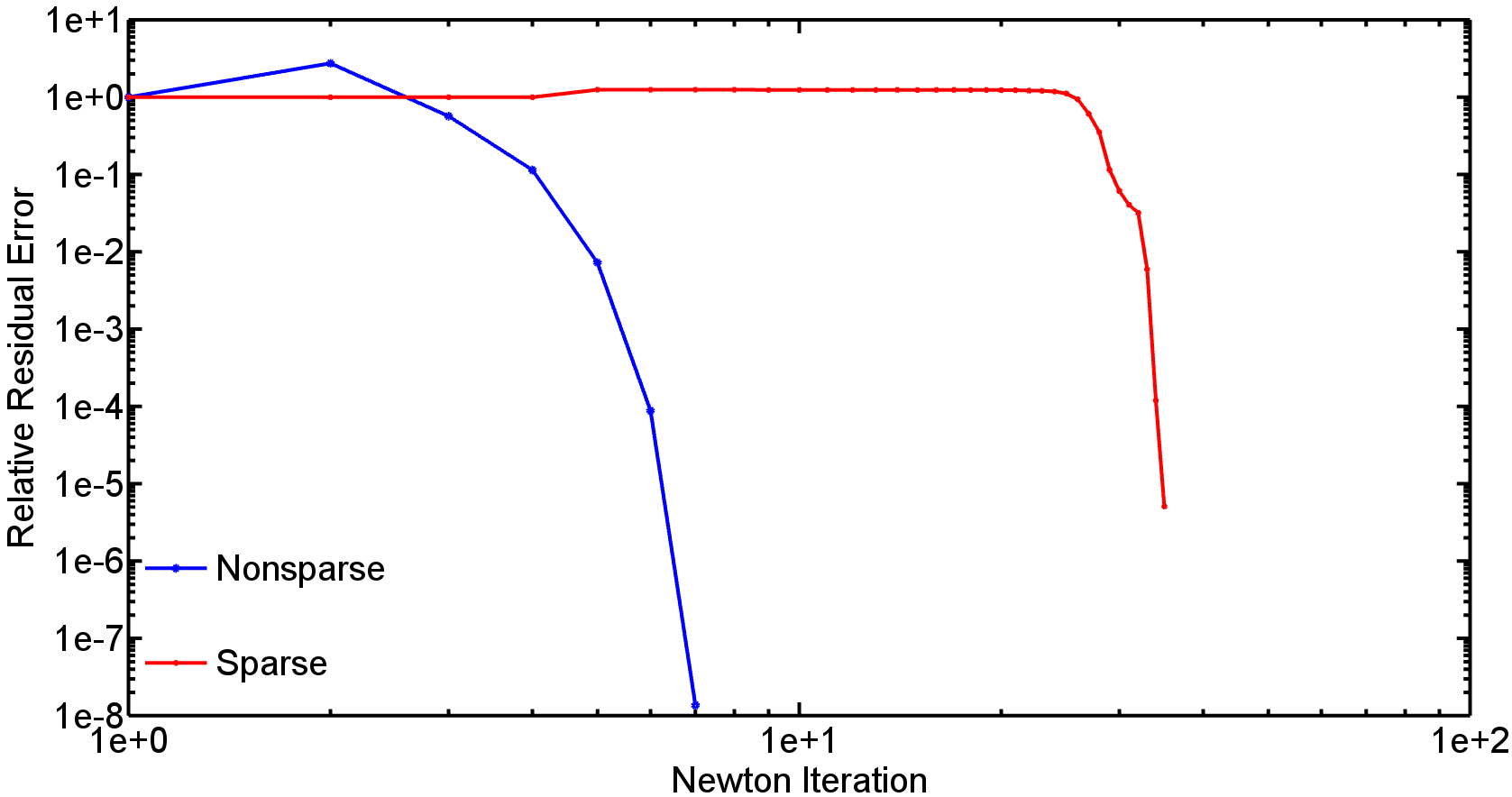}
\caption{Example 2, relative residual error reduction in the semismooth
Newton method for the nonsparse and sparse control in the turning wave
example.} 
\label{fig:wavenewton}
\end{figure}

To solve the nonlinear first order necessary optimality system, we needed $7$
and $35$ semismooth Newton iterations for the nonsparse and sparse optimal
control, respectively; see Fig.~\ref{fig:wavenewton}. The total system
assembling and solving time is about $3$ and $12.7$ hours on the desktop 
computer, respectively.

\begin{figure}[h]
  \centering
  \includegraphics[scale=0.25]{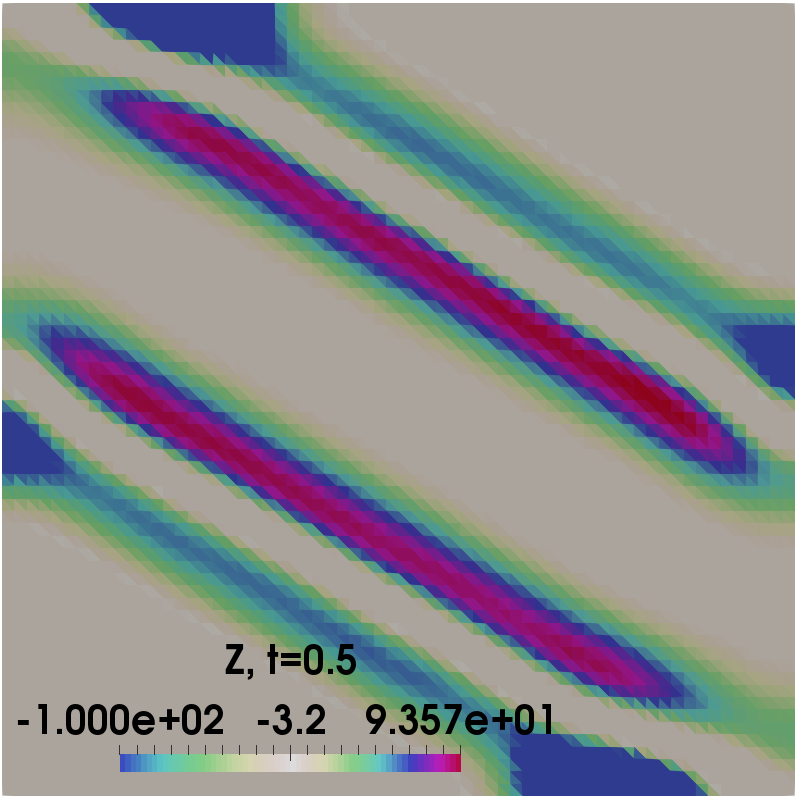}
  \includegraphics[scale=0.25]{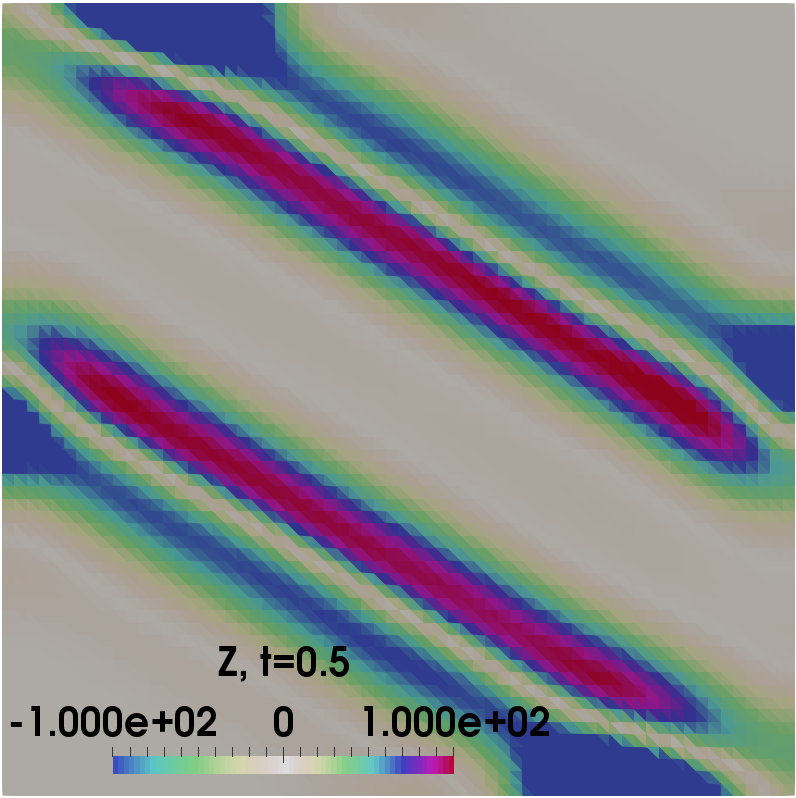}
  \includegraphics[scale=0.25]{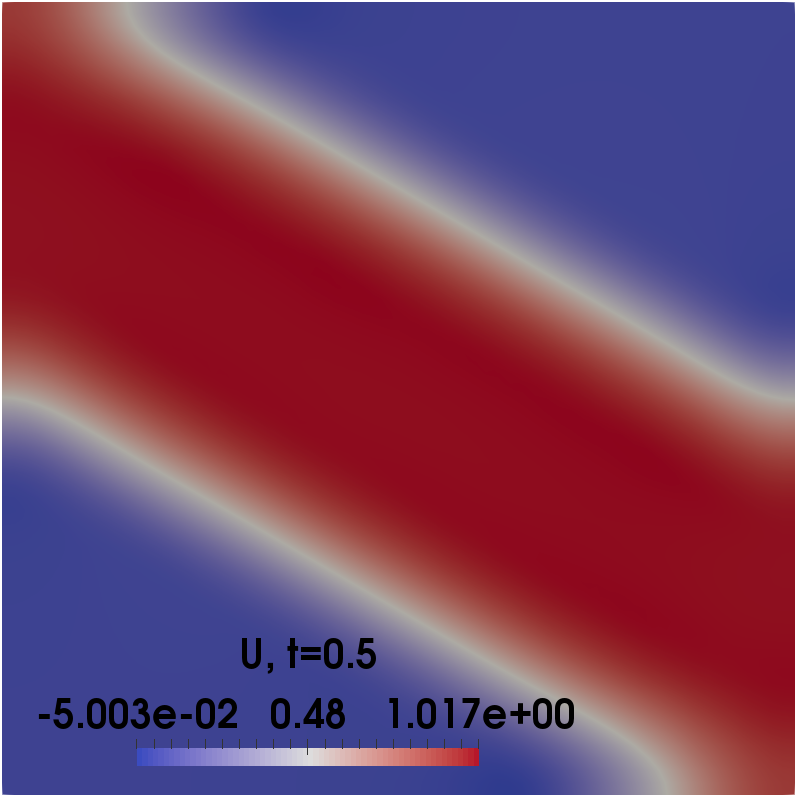}
  \includegraphics[scale=0.25]{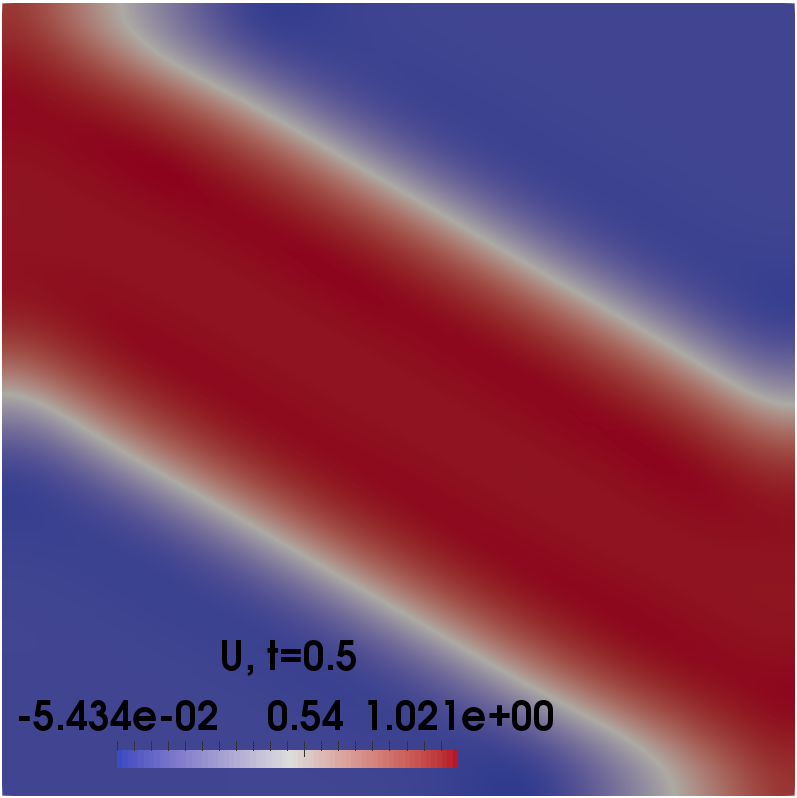}
  \caption{Example 2, plots of the sparse (left) and nonsparse (right)
    controls (in the first row) and the associated states (in the second row)
    at time $t=0.5$ for the turning wave example.} \label{fig:wavesparsecontrol}
\end{figure}

\begin{figure}[h]
  \centering
  \includegraphics[scale=0.25]{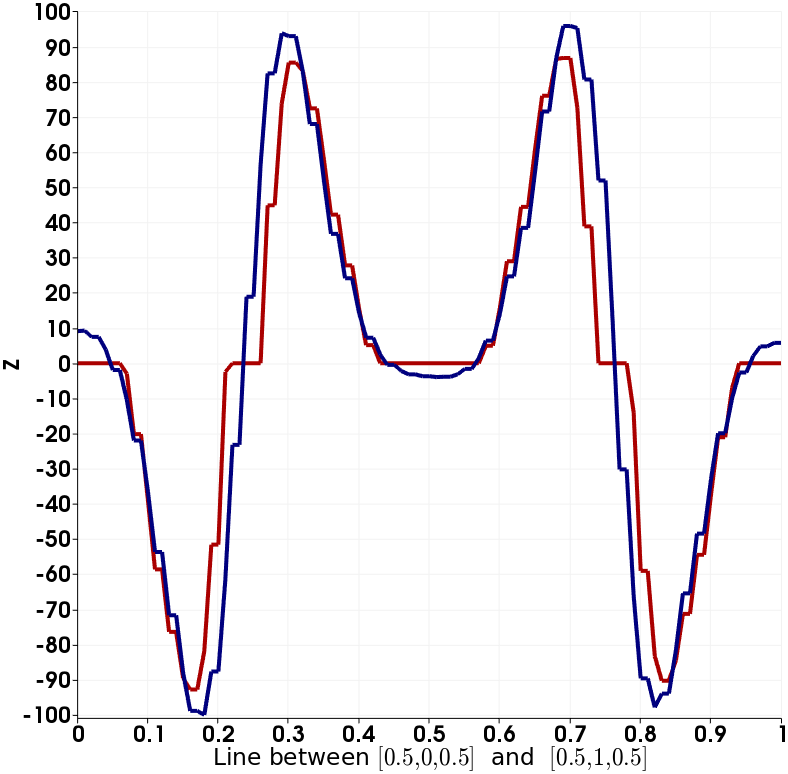}
  \includegraphics[scale=0.25]{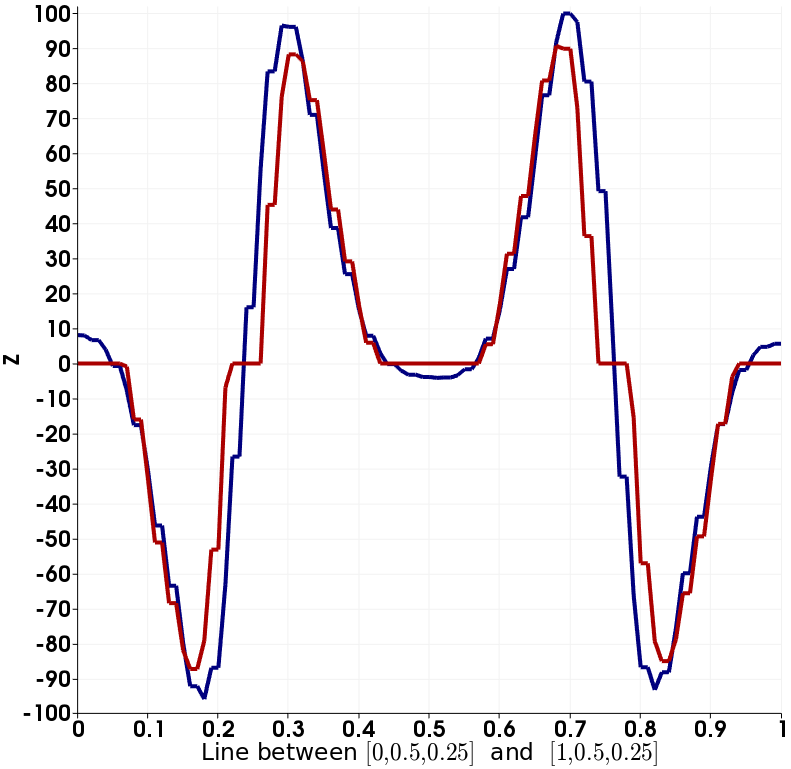}
  \includegraphics[scale=0.25]{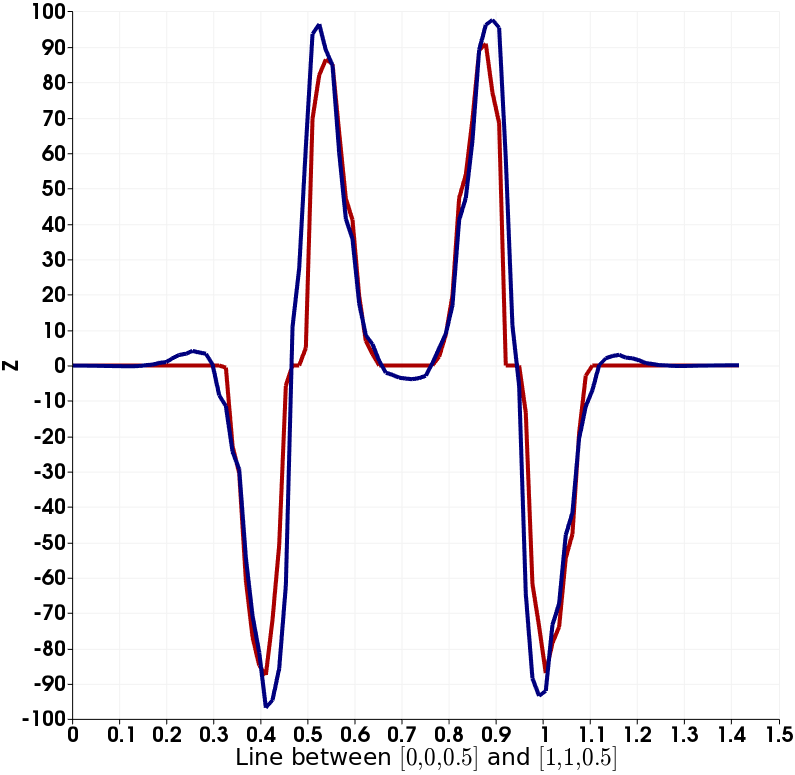}
  \caption{Example 2, comparisons of sparse (red) and nonsparse (blue)
    controls along the line between $[0.5, 0, 0.5]$ and $[0.5, 1, 0.5]$,
    between $[0, 0.5, 0.25]$ and $[1, 0.5, 0.25]$, and between $[0, 0, 0.5]$
    and $[1, 1, 0.5]$ for the turning wave example.} 
\label{fig:wavesparsecontrolline}  
\end{figure}

The numerical solutions of sparse and nonsparse controls as well as
associated optimal states are illustrated in
Fig.~\ref{fig:wavesparsecontrol}. We clearly see certain sparsity of our 
optimal sparse control as compared to pure $L^2$-regularization, without too
much precision loss of the associated state to the target. A closer look to the
sparse and nonsparse controls confirms that our optimal sparse
control exhibits sparsity with respect to the spatial direction; see
Fig.~\ref{fig:wavesparsecontrolline}.   

\begin{figure}[h]
  \centering
  \includegraphics[scale=0.248]{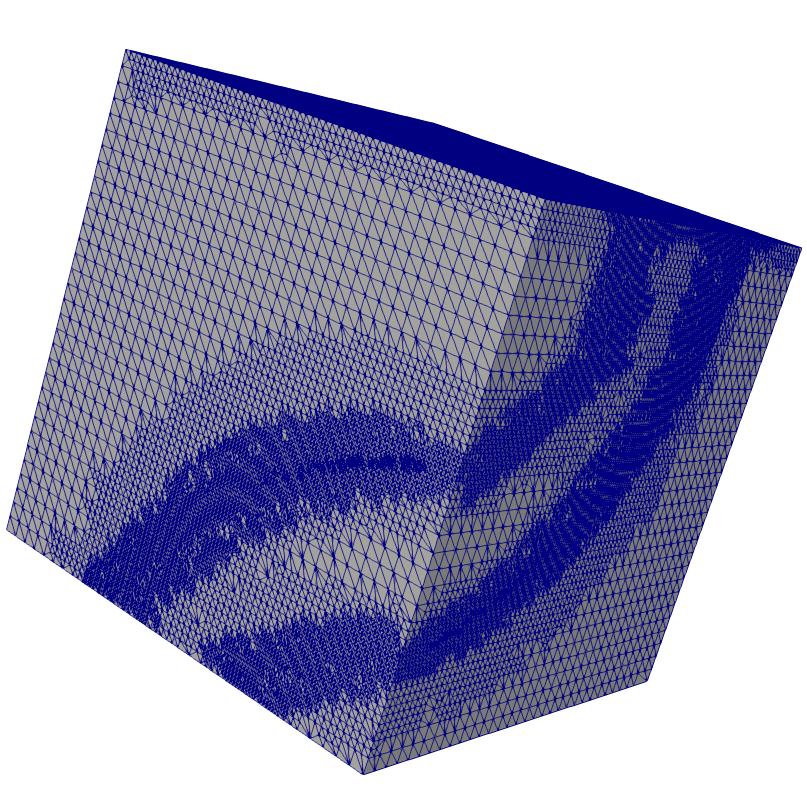}
  \includegraphics[scale=0.25]{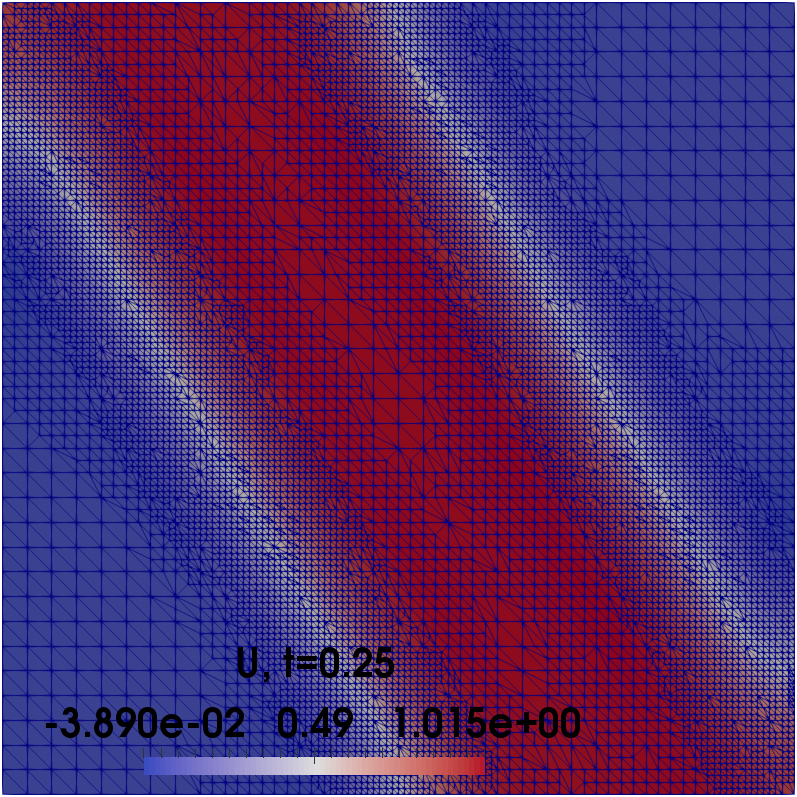}
  \includegraphics[scale=0.25]{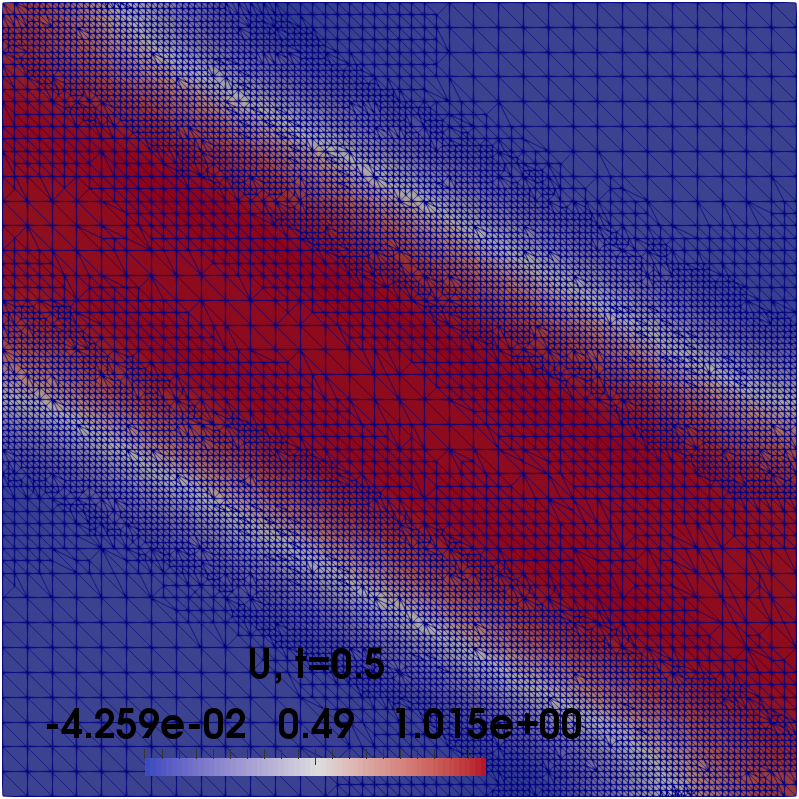}
  \includegraphics[scale=0.25]{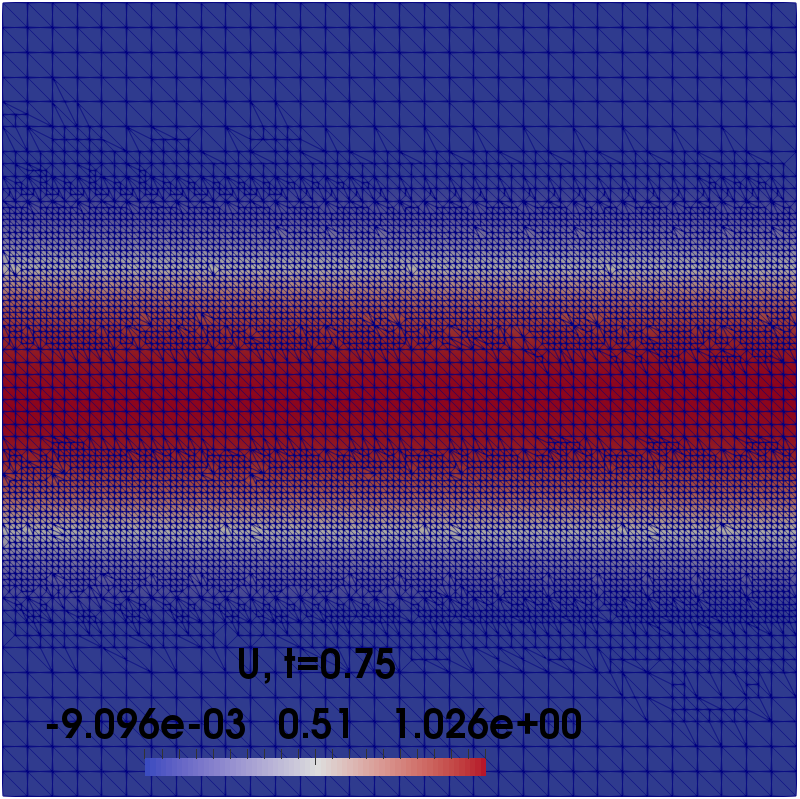}
  \caption{Example 2, plots of the adaptive space-time mesh (top-left) at the
    $6$th step, and the meshes on the cutting plans for times $t=0.25$, $0.5$, and $0.75$.} \label{fig:adaptmeshsparsecontrol}
\end{figure}

Instead of an uniform refinement, we may also adopt an adaptive strategy. In this
way, local refinements are made in the region where the
solution shows a more local character, whereas coarser meshes appear in the 
other region. For
example, in the optimal sparse control case ($\varrho=10^{-6}$, $\mu=10^{-4}$), we 
start from an initial mesh with $729$ grid points, $9$ in each spatial and the
temporal direction. We use a residual based error indicator for the coupled
state and adjoint state system to guide our adaptive mesh refinement, similar
to the approach \cite{SteinbachYang}. After the $6$th adaptive octasection
refinement \cite{Jurgen}, the mesh contains $1,053,443$ grid points; see the
adaptive space-time mesh and the meshes on the cutting plans at different times
in Fig.~\ref{fig:adaptmeshsparsecontrol}. As we observe, the adaptive 
refinements follow the rotation of the wave front of the state.   

\section{Conclusions}\label{sec:con}
In this work, we have considered a space-time Petrov-Galerkin finite element
method on fully unstructured simplicial meshes for semilinear parabolic 
optimal sparse control problems. The objective functional involves the 
well-known $L^1$-norm of the control in addition to
the standard $L^2$-regularization term. The proposed
method is able to capture spatio-temporal sparsity, which has been confirmed
by our numerical experiments. 
A rigorous convergence and error analysis of our
space-time Petrov-Galerkin finite element methods for such optimal sparse
control problems is left for future work.     

\bibliographystyle{plain}
\bibliography{LangerSteinbachTroeltzschYang2020}

\end{document}